\documentclass[12pt]{article}
\usepackage{amssymb}
\usepackage{amsmath}
\usepackage{array}
\usepackage[russian]{babel}
\usepackage[pdftex]{graphics}
\usepackage[X2,T2A]{fontenc}
\usepackage[cp1251]{inputenc}
\usepackage{mathrsfs}
\usepackage{longtable}
\usepackage{mathrsfs}
\usepackage{textcase}

\DeclareMathOperator{\mmod}{mod}

\newcommand{\prsum}{\mathop{{\sum}'}}

\multlinegap=0em \DeclareFontFamily{T1}{msb}{}
\DeclareFontShape{T1}{msb}{m}{ol}{<5> <6> <7> <8> <9> gen * msbm
<10> <10.95> <12> <14.4> <17.28> <20.74> <24.88> msbm10}{}
\DeclareSymbolFont{AMSb}{T1}{msb}{m}{ol} \multlinegap=0em

\renewcommand{\S}{\mathhexbox278}

\begin{document}

\begin{center}
{\rmfamily\bfseries\normalsize An upper bound for the second moment of the length of the period \\ of the continued fraction expansion for $\boldsymbol{\sqrt{d}}$}
\end{center}

\begin{center}
{\normalsize M.A.~Korolev\footnote{Steklov Mathematical Institute of Russian Academy of Sciences. E-mail: \texttt{korolevma@mi-ras.ru}, \texttt{hardy\_ramanujan@mail.ru}}}
\end{center}

\vspace{1cm}

\textbf{\S 1. Introduction.}
\vspace{0.5cm}

Suppose that $d\geqslant 1$ is an integer. If $d$ is not a perfect square, we define $T(d)$ as the length of the minimal period of the simple continued fraction expansion for $\sqrt{d}$. Otherwise, we put $T(d) = 0$. It is known (see \cite{Hickerson_1973}) that $T(d)\leqslant g(d)$, where $g(d)$ denotes the number of pairs of integers $m,q\geqslant 1$ such that
\[
m<\sqrt{d}, \quad |q-\sqrt{d}|<m,\quad q\,|\, (d-m^{2}).
\]
In the recent paper \cite{Battistoni_Grenie_Molteni_2024}, F.~Battistoni, L.~Greni\'{e} and G.~Molteni established the following relations for the first and second moments of the function
$g(d)$:
\begin{align}
& \sum\limits_{d\leqslant x}g(d) = \frac{4}{3}\,(\ln{2})\,x^{3/2} -2x-2\sqrt{x}+\theta(x+4\sqrt{x}),\quad 0\leqslant\theta\leqslant 1, \label{lab_01} \\
& \sum\limits_{d\leqslant x}g^{2}(d)\leqslant 11.9x^{2} + 5x^{3/2}\ln^{2}(4e^{4}x), \label{lab_02}
\end{align}
These formulas hold true for any $x>1$. The bound (\ref{lab_02}) together with (\ref{lab_01}) and Cauchy inequality imply the estimates
\begin{align*}
& \sum\limits_{d\leqslant x}T^{2}(d)\leqslant 11.9x^{2} + 5x^{3/2}\ln^{2}(4e^{4}x), \\
& \sum\limits_{x<d\leqslant 2x}T^{2}(d)\leqslant 47x^{2} + O\bigl(x^{3/2}\ln^{2}{x}\bigr).
\end{align*}
The last one was used in \cite{Battistoni_Grenie_Molteni_2024} for the proof of the inequality
\begin{equation}\label{lab_03}
\#D(x;\alpha)\leqslant \frac{47+o(1)}{\alpha^{2}}\,x
\end{equation}
for the cardinality of the set $D(x;\alpha)$ of $d$, $x<d\leqslant 2x$ satisfying the condition $T(d)>\alpha\sqrt{x}$. The inequality (\ref{lab_03}) improves the bound
\begin{equation*}
\#D(x;\alpha)\leqslant \frac{c+o(1)}{\ln^{2\mathstrut}{\alpha}}\,x,\quad c>0,
\end{equation*}
which belongs to A.M.~Rockett and P. Sz\"{u}sz \cite{Rockett_Szusz_1990}.

The authors of \cite{Battistoni_Grenie_Molteni_2024} express the second moment of $g(d)$ in the form
\[
W = \sum\limits_{\substack{1\leqslant m_{i}<\sqrt{x} \\ i = 1,2}}\;\;\sum\limits_{\substack{1\leqslant q_{i}<\sqrt{x} + m_{i} \\ i = 1,2}}
\;\;\sum\limits_{\substack{k_{1}, k_{2}\,:\,|k_{i}-q_{i}|<2m_{i},\, i=1,2 \\ m_{1}^{2}+k_{1}q_{1} = m_{2}^{2}+k_{2}q_{2} \leqslant x}}1.
\]
Then, omitting some restrictions on variables $k_{i}$, they arrived to the estimate (\ref{lab_02}).
In the end of \cite{Battistoni_Grenie_Molteni_2024}, the authors show briefly how to replace (\ref{lab_02}) by more precise bound with $8x^{2}+O\bigl(x^{3/2}\ln^{4}{x}\bigr)$ in the right-hand side.

The purpose of the present paper is to find an asymptotic formula for the sum $W$. Our main assertion is the following
\vspace{0.5cm}

\textsc{Theorem.} \emph{Let $0<\varepsilon<\tfrac{1}{54}$ be any fixed constant. Then, for $x\to +\infty$, the following asymptotic holds:}
\[
W = c_{0}x^{2} + O\bigl(x^{2-1/54+\varepsilon}\bigr),\quad c_{0} = \frac{13}{14}\frac{\zeta(2)}{\zeta(3)}\,\mathfrak{S},
\]
\emph{where}
\begin{multline*}
\mathfrak{S}  = \int_{0}^{1}du_{1}\int_{0}^{1}du_{2}\int_{0}^{1+u_{1}}dv_{1}\int_{0}^{1+u_{2}}\varphi(u_{1},u_{2},v_{1},v_{2})dv_{2},\quad \varphi(u_{1},u_{2},v_{1},v_{2}) =\\
= \frac{1}{v_{1}v_{2}}\max{\bigl(0,\min{(1,(u_{1}+v_{1})^{2},(u_{2}+v_{2})^{2})} - \max{(u_{1}^{2},u_{2}^{2},(u_{1}-v_{1})^{2},(u_{2}-v_{2})^{2})}\bigr)}
\end{multline*}

\textsc{Remark.} The approximate calculations show that $\mathfrak{S} = 0.959\ldots, c_{0} = 1.218\ldots$. After the submission of the 1st version of this paper, 11.05.2024, the authors of 
\cite{Battistoni_Grenie_Molteni_2024} kindly informed me that they had found the closed expression for the constant $\mathfrak{S}$. It turned out that $\mathfrak{S} = 2(\ln{2})^{2} = 0.9609060278\ldots$.
\vspace{0.5cm}

\textsc{Corollary.} \emph{Under the assumptions of Theorem, the following inequality holds for any $\alpha>0$:}
\[
\#D(x;\alpha)\leqslant \frac{c_{1}x}{\alpha^{2}}\bigl(1 + O(x^{-1/54+\varepsilon})\bigr),\quad c_{1} = 3c_{0} = 3.654\ldots,
\]
\vspace{0.3cm}

\textbf{\S 2. Auxiliary assertions.}
\vspace{0.3cm}

In this section, we give some auxiliary lemmas which are necessary for the following.
\vspace{0.3cm}

\textsc{Lemma 1.} \emph{Suppose that the function $f$ belongs to the class $C^{1}[a,b]$ and the number of the segments of its monotonicity does not exceed $k$} ($k\geqslant 1$). \emph{Then}
\[
\sum\limits_{a<n\leqslant b}f(n) = \int_{a}^{b}f(x)dx + \theta(k+1)\max_{a\leqslant x\leqslant b}{|f(x)|},\quad |\theta|\leqslant 1.
\]
This assertion easily follows from Euler summation formula.

Let $\varrho(u) = 1/2 - \{u\}$. Then for any $H>1$ one has $\varrho(u) = \varrho_{H}(u)+r_{H}(u)$, where
\begin{multline*}
\varrho_{H}(u) = \sum\limits_{0<|h|\leqslant H}\frac{e(hu)}{2\pi ih},\quad |r_{H}(u)|\ll \min{\Bigl(1, \frac{1}{H\|u\|}\Bigr)},\\
\|u\| = \min{(\{u\},1-\{u\})},\quad e(v) = e^{2\pi iv},
\end{multline*}
and the implied constant is absolute (see, for example, \cite[\S 5]{Hooley_1963}). To handle with the <<remainders>> $r_{H}(u)$, we need the following assertion.
\vspace{0.3cm}

\textsc{Lemma 2.} \emph{Let $H>3$. Then there exists a $1$-periodic function $\varrho_{H}(u)$ with the following properties:}

1) $\min{\bigl(1, (H\|u\|)^{-1}\bigr)} \leqslant 2\vartheta_{H}(u)$;

2) \emph{Fourier-series expansion of $\vartheta_{H}(u)$ has the form}
\[
\sum\limits_{h=-\infty}^{+\infty}C(h)e(hu),\quad \textit{where }\quad |C(h)|\leqslant \frac{2\ln{H}+1}{H}\quad\textit{for any}\quad h;
\]

3) \emph{for any fixed $A>3$ and $0<\varepsilon < 0.5$ one has}
\[
\sum\limits_{|h|>H^{1+\varepsilon}}|C(h)|\ll_{A,\varepsilon} H^{-A}.
\]

For the proof, see the paper \cite{Tolev_2005} of D.I.~Tolev.
\vspace{0.3cm}

\textsc{Lemma 3.} \emph{Denote by $\rho(\Delta)$ the number of solutions of the congruence}
\[
x^{2}\equiv y^{2}\pmod{\Delta},\quad 1\leqslant x,y \leqslant \Delta.
\]
\emph{Then $\rho(\Delta)$ is a multiplicative function, and the following relations hold:}
\begin{align*}
& \text{(a)} \quad \rho(\Delta) = \Delta\sum\limits_{\delta | \Delta}\frac{\varepsilon(\delta)}{\delta}\,\varphi(\delta),\quad \varepsilon(\delta) =
\begin{cases}
1, & \textit{if } \delta\equiv 1\pmod{2},\\
2, & \textit{if } \delta\equiv 0\pmod{4},\\
0, & \textit{if } \delta\equiv 2\pmod{4};
\end{cases} \\
& \text{(b)} \quad \sum\limits_{\Delta = 1}^{+\infty}\frac{\rho(\Delta)}{\Delta^{3}} = \frac{13}{14}\,\frac{\zeta^{2}(2)}{\zeta(3)}.
\end{align*}
\emph{Moreover, the following inequality holds: $\rho(\Delta)\leqslant 2\Delta\tau(\Delta)$, where $\tau(\Delta)$ is the divisor function.}

\textsc{Proof.} In fact, this is Lemma 2 from \cite{Battistoni_Grenie_Molteni_2024}. Here we give an alternative proof. Obviously,
\[
\rho(\Delta) = \sum\limits_{x,y = 1}^{\Delta}\frac{1}{\Delta}\sum\limits_{c=1}^{\Delta}e\Bigl(\frac{c}{\Delta}(x^{2}-y^{2})\Bigr) = \frac{1}{\Delta}\sum\limits_{c=1}^{\Delta}|S(\Delta,c)|^{2},\quad S(q,a) = \sum\limits_{x=1}^{q}e\Bigl(\frac{ax^{2}}{q}\Bigr).
\]
Let $(c,\Delta) = \omega$ and $\Delta/\omega = \delta, c/\omega = e$. Then
\[
S(\Delta,c) = \sum\limits_{x=1}^{\Delta}e\Bigl(\frac{ex^{2}}{\delta}\Bigr) = \frac{\Delta}{\delta}\sum\limits_{x=1}^{\delta}e\Bigl(\frac{ex^{2}}{\delta}\Bigr) = \frac{\Delta}{\delta}S(\delta,e).
\]
Therefore,
\[
\rho(\Delta)  = \frac{1}{\Delta}\sum\limits_{\omega\,|\,\Delta}\sum\limits_{\substack{1\leqslant c\leqslant \Delta \\ (c,\Delta) = \omega}}\frac{\Delta^{2}}{\delta^{2}}|S(\delta,e)|^{2}
= \Delta\sum\limits_{\delta | \Delta} \frac{1}{\delta^{2}}\sum\limits_{\substack{1\leqslant e \leqslant \delta \\ (e,\delta)=1}}|S(\delta,e)|^{2}.
\]
Since $(e,\delta)=1$ then $|S(\delta,e)|^{2} = \varepsilon(\delta)\delta$ (see, for example, \cite[Ch. I, \S 3, Theorem 3]{Korobov_2013}), we obtain the formula (a) and the inequality $\rho(\Delta)\leqslant 2\Delta\tau(\Delta)$. Next, for prime $p\geqslant 3$ this formula yields
\begin{multline*}
\frac{\rho(p^{\alpha})}{p^{3\alpha}} = \frac{1}{p^{2\alpha}} + \biggl(1-\frac{1}{p}\biggr)\frac{\alpha}{p^{2\alpha}},\\
\sum\limits_{\alpha = 0}^{+\infty}\frac{\rho(p^{\alpha})}{p^{3\alpha}} = 1 + \frac{1}{p^{2\mathstrut}-1} + \frac{p(p-1)}{(p^{2\mathstrut}-1)^{2\mathstrut}} = \frac{p(p^{3}-1)}{(p^{2\mathstrut}-1)^{2\mathstrut}} = \frac{1-1/p^{3}}{(1-1/p^{2})^{2}}.
\end{multline*}
Finally,
\[
\rho(2^{\alpha}) = \alpha 2^{\alpha},\quad \frac{\rho(2^{\alpha})}{2^{3\alpha}} = \frac{13}{9} = \frac{13}{14}\cdot \frac{1-1/2^{3}}{(1-1/2^{2})^{2}}.
\]
Hence,
\[
\sum\limits_{\Delta = 1}^{+\infty}\frac{\rho(\Delta)}{\Delta^{3}} = \prod_{p}\sum\limits_{\alpha = 0}^{+\infty}\frac{\rho(p^{\alpha})}{p^{3\alpha}} = \frac{13}{14}\prod_{p}\frac{1-1/p^{3}}{(1-1/p^{2})^{2}} = \frac{13}{14}\,\frac{\zeta^{2}(2)}{\zeta(3)}.
\]

\textsc{Lemma 4.} \emph{Let $L,M,N>1$ and $L<L_{1}\leqslant 2L$, $M<M_{1}\leqslant 2M$, $N<N_{1}\leqslant 2N$ and let $\mathbf{a} = \{a_{m}\}$, $\mathbf{b} = \{b_{n}\}$, $\mathbf{c} = \{c_{\ell}\}$ be any sequences defined
for $M<m\leqslant M_{1}$, $N<n\leqslant N_{1}$ and $L<\ell\leqslant L_{1}$. Further, let $\vartheta\ne 0$ is a fixed real number. Finally, suppose that the function $f_{\ell,\vartheta}(x,y)$ belongs to the class $C^{1}(\mathbb{R}_{+}^{2})$ for any $\ell$ and satisfies to the conditions}
\begin{equation}\label{lab_04}
\frac{\partial_{\ell,\vartheta}f}{\partial x}(x,y)\ll \frac{X}{x^{2}y},\quad \frac{\partial_{\ell,\vartheta}f}{\partial y}(x,y)\ll \frac{X}{x y^{2}}
\end{equation}
\emph{for $M\leqslant x\leqslant M_{1}$, $N\leqslant y\leqslant N_{1}$. Then, for any fixed $\varepsilon$, the sum}
\[
S = \sum\limits_{L<\ell\leqslant L_{1}}\mathop{\sum_{M<m\leqslant M_{1}}\sum_{N<n\leqslant N_{1}}}\limits_{(m,n)=1}a_{m}b_{n}c_{\ell}\,e\biggl(\vartheta\,\frac{\ell\,\overline{m}}{n} + f_{\ell,\vartheta}(m,n)\biggr)
\]
\emph{satisfies the estimate}
\begin{multline*}
S\ll_{\varepsilon}(LMN)^{\varepsilon}\|\mathbf{a}\|\cdot \|\mathbf{b}\|\cdot \|\mathbf{c}\|\cdot T\Bigl((LMN)^{7/20}(M+N)^{1/4} + (LMN)^{3/8}(M+N)^{1/8}L^{1/8}\Bigr),\\
T = \biggl(1+\frac{|\vartheta|L+X}{MN}\biggr)^{1/2},\quad \|\mathbf{a}\| = \Bigl(\,\sum\limits_{M<m\leqslant M_{1}}|a_{m}|^{2} \Bigr)^{1/2}\quad (\textit{similarly for } \|\mathbf{b}\|, \|\mathbf{c}\|).
\end{multline*}
\renewcommand{\refname}{\normalsize{References}}

For the proof, see \cite[Theorem 1 and Remark 1]{Bettin_Chandee_2015}.
\pagebreak

\textbf{\S 3. Initial transformations of the sum $\boldsymbol{W}$.}
\vspace{0.5cm}

First we estimate the contribution coming to $W$ from the terms with $q_{1} = q_{2}$. If $q_{1} = q_{2} = q$ then the relation $m_{1}^{2}+k_{1}q = m_{2}^{2}+k_{2}q$ implies the congruence
$m_{1}^{2}\equiv m_{2}^{2}\pmod{q}$. Given such $q, m_{1}, m_{2}$, the value of $k_{2}$ is determined uniquely by $k_{1}$: $k_{2} = k_{1} + (m_{1}^{2}-m_{2}^{2})/q$. Hence, if we omit the condition $|k_{2}-q|<2m_{2}$ then we bound the contribution under considering by the sum
\[
R = \sum\limits_{q\leqslant 2\sqrt{x}}\sum\limits_{\substack{m_{i}\leqslant \sqrt{x},\; i=1,2 \\ m_{1}^{2}\equiv m_{2}^{2}\;(\mmod{q})}}\sum\limits_{q-2m_{1}<k_{1}< q+2m_{1}}1
\leqslant 4\sqrt{x}\sum\limits_{q\leqslant 2\sqrt{x}}\sum\limits_{\substack{m_{i}\leqslant \sqrt{x},\; i=1,2 \\ m_{1}^{2}\equiv m_{2}^{2}\;(\mmod{q})}}1.
\]
Splitting the sums over $m_{1}, m_{2}$ into arithmetic progressions modulo $q$ and using Lemma 3 we get
\begin{multline*}
R \leqslant 4\sqrt{x}\sum\limits_{q\leqslant 2\sqrt{x}}\sum\limits_{\substack{\xi_{1},\xi_{2} = 0 \\ \xi_{1}^{2}\equiv \xi_{2}^{2}\;(\mmod{q})}}^{q-1}
\sum\limits_{\substack{m_{i}\leqslant \sqrt{x} \\ m_{i}\equiv \xi_{i}\;(\mmod{q}) \\ i=1,2}}1
\leqslant 4\sqrt{x}\sum\limits_{q\leqslant 2\sqrt{x}}\sum\limits_{\substack{\xi_{1},\xi_{2} = 0 \\ \xi_{1}^{2}\equiv \xi_{2}^{2}\;(\mmod{q})}}^{q-1}\biggl(\frac{\sqrt{x}}{q}+1\biggr)^{2}\ll \\
\ll x^{3/2}\sum\limits_{q\leqslant 2\sqrt{x}}\frac{\rho(q)}{q^{2}} \ll x^{3/2}\sum\limits_{q\leqslant 2\sqrt{x}}\frac{\tau(q)}{q}\ll x^{3/2}(\ln{x})^{2}.
\end{multline*}
Hence, since the sum $W$ is symmetrical with respect to $q_{1}, q_{2}$, we find that
\[
W = 2W_{1} + O\bigl(x^{3/2}(\ln{x})^{2}\bigr),\quad
W_{1} =  \sum\limits_{\substack{1\leqslant m_{i}<\sqrt{x} \\ i = 1,2}}\;\;\sum\limits_{\substack{1\leqslant q_{i}<\sqrt{x} + m_{i} \\ i = 1,2 \\ q_{1}\leqslant q_{2}}}
\;\;\sum\limits_{\substack{k_{1}, k_{2}\,:\,|k_{i}-q_{i}|<2m_{i},\, i=1,2 \\ m_{1}^{2}+k_{1}q_{1} = m_{2}^{2}+k_{2}q_{2} \leqslant x}}1.
\]

Now we transform the sum $W_{1}$. The relation $m_{1}^{2} + k_{1}q_{1} = m_{2}^{2} + k_{2}q_{2}$ implies the congruence
\begin{equation}\label{lab_05}
m_{1}^{2} + k_{1}q_{1}\equiv m_{2}^{2}\pmod{q_{2}}.
\end{equation}
If the quadruple $m_{1}, m_{2}, q_{1}, k_{1}$ satisfies (\ref{lab_05}) then the number $k_{2} = (k_{1}q_{1}+m_{1}^{2}-m_{2}^{2})/q_{2}$ is integer. In this case, the condition
$|k_{2}-q_{2}|<2m_{2}$ (or, that is the same, $q_{2}-2m_{2}<k_{2}<q_{2}+2m_{2}$) takes the form
\[
q_{2}-2m_{2}<\frac{k_{1}q_{1}+m_{1}^{2}-m_{2}^{2}}{q_{2}} <q_{2}+2m_{2}
\]
or, equivalently, the form
\begin{equation}\label{lab_06}
\frac{(q_{2}-m_{2})^{2}-m_{1}^{2}}{q_{1}}<k_{1}<\frac{(q_{2}+m_{2})^{2}-m_{1}^{2}}{q_{1}}.
\end{equation}
Further, the inequalities $m_{1}^{2}+k_{1}q_{1}\leqslant x$, $k_{1}>0$ imply
\begin{equation}\label{lab_07}
0<k_{1}\leqslant \frac{x-m_{1}^{2}}{q_{1}}.
\end{equation}
Next, the condition $k_{2}>0$ takes the form
\begin{equation}\label{lab_08}
k_{1} > \frac{m_{2}^{2}-m_{1}^{2}}{q_{1}}.
\end{equation}
Finally, we write the inequality $|k_{1}-q_{1}|<2m_{1}$ in the form
\begin{equation}\label{lab_09}
q_{1}-2m_{1}<k_{1}<q_{1}+2m_{1}.
\end{equation}
Now we replace the upper bounds for $k_{1}$ in (\ref{lab_06}) and (\ref{lab_09}) by weak inequalities and estimate the contribution coming to $W_{1}$ of such replacement.

Let $k_{1} = q_{1}+2m_{1}$. In view of (\ref{lab_05}), we have
\begin{equation}\label{lab_10}
(q_{1}+m_{1})^{2}-m_{2}^{2}\equiv 0\pmod{q_{2}}\quad \text{and therefore }\quad q_{2}\,|\,(q_{1}+m_{1})^{2}-m_{2}^{2}.
\end{equation}
The contribution coming from $q_{2}+m_{1}\neq m_{2}$ is less than
\[
\sum\limits_{m_{1},m_{2}\leqslant \sqrt{x}}\;\sum\limits_{\substack{q_{1}\leqslant 2\sqrt{x} \\ q_{2}\neq m_{2}-m_{1}}}\tau\bigl(|(q_{1}+m_{1})^{2}-m_{2}^{2}|\bigr)\ll x^{3/2+\varepsilon}.
\]
The terms with $q_{1} = m_{2}-m_{1}>0$ yield
\[
\ll \sum\limits_{m_{1},m_{2}\leqslant \sqrt{x}}\;\sum\limits_{q_{2}\leqslant 2\sqrt{x}}1\ll x^{3/2}.
\]
Similarly, suppose that $k_{1} = ((q_{2}+m_{2})^{2}-m_{1}^{2})/q_{1}$ is integer. Then (\ref{lab_05}) holds automatically. We may also suppose that $q_{2}\neq m_{1} - m_{2}$: otherwise, we get $k_{1} = 0$. Then $q_{2}$ divides $(q_{2}+m_{2})^{2}-m_{1}^{2}\neq 0$, and the corresponding contribution to $W_{1}$ is less than
\[
\sum\limits_{m_{1},m_{2}\leqslant \sqrt{x}}\; \sum\limits_{\substack{q_{2}\leqslant 2\sqrt{x} \\ q_{2}\neq m_{1}-m_{2}}}\tau\bigl(|(q_{2}+m_{2})^{2}-m_{1}^{2}|\bigr)\ll x^{3/2+\varepsilon}.
\]
Thus we get $W_{1} = W_{2} + O\bigl(x^{3/2+\varepsilon}\bigr)$, where
\begin{multline*}
W_{2} = \prsum\limits_{\substack{m_{i}\leqslant \sqrt{x} \\ i = 1,2}}\;\prsum\limits_{\substack{1\leqslant q_{i}\leqslant \sqrt{x}+m_{i}, \; i=1,2 \\ q_{1}\leqslant q_{2}}}\;\sum\limits_{\substack{f<k\leqslant g \\ kq_{1}\equiv m_{2}^{2}-m_{1}^{2}\pmod{q_{2}}}}1,\\
f = f(x;m_{1},m_{2},q_{1},q_{2}) = \max{\biggl(0,\frac{m_{2}^{2}-m_{1}^{2}}{q_{1}},q_{1}-2m_{1},\frac{(q_{2}-m_{2})^{2}-m_{1}^{2}}{q_{1}}\biggr)},\\
g = g(x;m_{1},m_{2},q_{1},q_{2}) = \min{\biggl(\frac{x-m_{1}^{2}}{q_{1}},q_{1}+2m_{1},\frac{(q_{2}+m_{2})^{2}-m_{1}^{2}}{q_{1}}\biggr)}
\end{multline*}
and the prime sign means the summation over the set $\mathcal{A}(x)$ of quadruples $(m_{1},m_{2},q_{1},q_{2})$ of positive integers that satisfy the restrictions
$1\leqslant m_{i}\leqslant \sqrt{x}$, $1\leqslant q_{i}\leqslant \sqrt{x}+m_{i}$, $i=1,2$, $q_{1}\leqslant q_{2}$, and, moreover, obey the inequality
\begin{equation}\label{lab_11}
f(x;m_{1},m_{2},q_{1},q_{2}) < g(x;m_{1},m_{2},q_{1},q_{2}).
\end{equation}
Obviously, the functions $f$ and $g$ can be written in the form
\begin{multline*}
f = \frac{1}{q_{1}}\bigl(\max{\bigl\{m_{1}^{2},m_{2}^{2},(q_{1}-m_{1})^{2},(q_{2}-m_{2})^{2}\bigr\}}-m_{1}^{2}\bigr),\\
g = \frac{1}{q_{1}}\bigl(\min{\bigl\{x,(q_{1}+m_{1})^{2},(q_{2}+m_{2})^{2}\bigr\}}-m_{1}^{2}\bigr).
\end{multline*}
This means that the condition (\ref{lab_11}) is equivalent to the inequality
\[
\max{(f_{1},f_{2},f_{3},f_{4})}<\min{(g_{1},g_{2},g_{3})},
\]
where
\begin{multline*}
f_{1} = m_{1}^{2},\quad f_{2} = m_{2}^{2},\quad f_{3} = (q_{1}-m_{1})^{2},\quad f_{4} = (q_{2}-m_{2})^{2},\\
g_{1} = x,\quad g_{2} = (q_{1}+m_{1})^{2}, \quad g_{3} = (q_{2}+m_{2})^{2}.
\end{multline*}
If the tuple $(m_{1},m_{2},q_{1},q_{2})$ belongs to $\mathcal{A}(x)$ then $m_{1}, m_{2}$ satisfy the condition $m_{1}^{2}\equiv m_{2}^{2}\pmod{\Delta}$, where $\Delta = (q_{1},q_{2})$. Setting $q_{i} = \Delta \kappa_{i}$, we get $1\leqslant \kappa_{i}\leqslant (\sqrt{x}+m_{i})/\Delta$, $i = 1,2$, and $(\kappa_{1}, \kappa_{2})=1$. Therefore, the congruence (\ref{lab_05}) takes the form $k\Delta \kappa_{1}\equiv \Delta\nu \pmod{\Delta\kappa_{2}}$, and has solutions
\[
k\equiv \nu\bar{\kappa}_{1}\pmod{\kappa_{2}},\quad \text{where}\quad \nu = \frac{m_{2}^{2}-m_{1}^{2}}{\Delta},\quad \kappa_{1}\bar{\kappa}_{1}\equiv 1\pmod{\kappa_{2}}.
\]
Choosing $\varkappa$ from the relations $\varkappa\equiv \nu\bar{\kappa}_{1}\pmod{\kappa_{2}}$, $0\leqslant \varkappa<\kappa_{2}$ and setting $k = \tau\kappa_{2}+\varkappa$, we get
$(f-\varkappa)/\kappa_{2}<\tau\leqslant (g-\varkappa)/\kappa_{2}$. Thus the inner sum in $W_{2}$ is equal to
\[
\biggl[\frac{g-\varkappa}{\kappa_{2}}\biggr]\,-\,\biggl[\frac{f-\varkappa}{\kappa_{2}}\biggr].
\]
Hence,
\[
W_{2} = \sum\limits_{\Delta\leqslant 2\sqrt{x}}\;\prsum\limits_{\substack{1\leqslant m_{i}\leqslant\sqrt{x},\;i = 1,2 \\ m_{1}^{2}\equiv m_{2}^{2}\;(\mmod{\Delta})}}\;\prsum\limits_{\substack{1\leqslant \kappa_{i}\leqslant (\sqrt{x}+m_{i})/\Delta,\;i=1,2 \\ (\kappa_{1},\kappa_{2})=1,\; \kappa_{1}\leqslant \kappa_{2}}}\biggl(\biggl[\frac{g-\varkappa}{\kappa_{2}}\biggr]\,-\,\biggl[\frac{f-\varkappa}{\kappa_{2}}\biggr]\biggr).
\]
Further, we split the domains of $m_{i}$, $i=1,2$, to arithmetic progressions modulo $\Delta$ as follows:
\[
W_{2} = \sum\limits_{\Delta\leqslant 2\sqrt{x}}\;\sum\limits_{\substack{0\leqslant\xi_{1},\xi_{2}< \Delta \\ \xi_{1}^{2}\equiv \xi_{2}^{2}\;(\mmod{\Delta})}}\prsum\limits_{\substack{1\leqslant m_{i}\leqslant\sqrt{x},\\ m_{i}\equiv \xi_{i}\;(\mmod{\Delta}),\;i = 1,2}}\;\prsum\limits_{\substack{1\leqslant \kappa_{i}\leqslant (\sqrt{x}+m_{i})/\Delta,\;i=1,2 \\ (\kappa_{1},\kappa_{2})=1,\; \kappa_{1}\leqslant \kappa_{2}}}\biggl(\biggl[\frac{g-\varkappa}{\kappa_{2}}\biggr]\,-\,\biggl[\frac{f-\varkappa}{\kappa_{2}}\biggr]\biggr).
\]
Suppose now that $(\ln{x})^{4}\ll D<\sqrt[4\;]{x}$ (the precise value for $D$ will be chosen later) and estimate the contribution coming to $W_{2}$ from the terms with $D<\Delta\leqslant 2\sqrt{x}$. In the case $g>f$ we obviously have
\[
\frac{g-f}{\kappa_{2}}\leqslant \frac{g}{\kappa_{2}}\leqslant \frac{x-m_{1}^{2}}{q_{1}\kappa_{2}} < \frac{x}{\Delta\kappa_{1}\kappa_{2}}.
\]
By Lemma 3, such a contribution does not exceed
\begin{multline*}
x\sum\limits_{D<\Delta\leqslant 2\sqrt{x}}\frac{1}{\Delta}\sum\limits_{\substack{0\leqslant\xi_{1},\xi_{2}< \Delta \\ \xi_{1}^{2}\equiv \xi_{2}^{2}\;(\mmod{\Delta})}}\prsum\limits_{\substack{1\leqslant m_{i}\leqslant\sqrt{x},\\ m_{i}\equiv \xi_{i}\;(\mmod{\Delta}),\;i = 1,2}}\;\prsum\limits_{\substack{1\leqslant \kappa_{i}\leqslant 2\sqrt{x}/\Delta \\ i=1,2}}\frac{1}{\kappa_{1}\kappa_{2}} \ll x(\ln{x})^{2}\sum\limits_{D<\Delta\leqslant 2\sqrt{x}}\frac{\rho(\Delta)}{\Delta}\biggl(\frac{\sqrt{x}}{\Delta}\biggr)^{2} \\
\ll\,x^{2}(\ln{x})^{2}\sum\limits_{\Delta>D}\frac{\rho(\Delta)}{\Delta^{3}}\,\ll x^{2}(\ln{x})^{2}\sum\limits_{\Delta>D}\frac{\tau(\Delta)}{\Delta^{2}}\,\ll \frac{x^{2}(\ln{x})^{3}}{D}.
\end{multline*}
Next, suppose that $\ln{x}\ll C \leqslant \sqrt{x}/(2D)$ (the value of $C$ will also be chosen later) and estimate the contribution coming to $W_{2}$ from the terms with $\kappa_{1}\leqslant \sqrt{x}/(\Delta C)$.

In the case $g>f$ we have
\[
\frac{g-f}{\kappa_{2}}\leqslant\frac{g}{\kappa_{2}}\leqslant \frac{q_{1}+2m_{1}}{\kappa_{2}}\leqslant\frac{\sqrt{x}+3m_{1}}{\kappa_{2}} \leqslant \frac{4\sqrt{x}}{\kappa_{2}}.
\]
Therefore, such contribution does not exceed
\begin{multline*}
\sum\limits_{\Delta\leqslant D}\sum\limits_{\substack{0\leqslant\xi_{1},\xi_{2}< \Delta \\ \xi_{1}^{2}\equiv \xi_{2}^{2}\;(\mmod{\Delta})}}\sum\limits_{\substack{1\leqslant m_{i}\leqslant\sqrt{x},\\ m_{i}\equiv \xi_{i}\;(\mmod{\Delta}),\;i = 1,2}}\;\sum\limits_{1\leqslant \kappa_{1}\leqslant 2\sqrt{x}/(C\Delta)}
\sum\limits_{1\leqslant \kappa_{2}\leqslant 2\sqrt{x}/\Delta}\frac{4\sqrt{x}}{\kappa_{2}}\ll\\
\ll\frac{x\ln{x}}{C}\sum\limits_{\Delta\leqslant D}\frac{1}{\Delta}\sum\limits_{\substack{0<\xi_{1},\xi_{2}\leqslant \Delta \\ \xi_{1}^{2}\equiv \xi_{2}^{2}\;(\mmod{\Delta})}}
\sum\limits_{\substack{1\leqslant m_{i}\leqslant\sqrt{x},\\ m_{i}\equiv \xi_{i}\;(\mmod{\Delta}),\;i = 1,2}}1\ll  \frac{x\ln{x}}{C}\sum\limits_{\Delta\leqslant D}\frac{\rho(\Delta)}{\Delta}\biggl(\frac{\sqrt{x}}{\Delta}\biggr)^{2}\ll \frac{x^{2}\ln{x}}{C}.
\end{multline*}
Finally, since $[w]-[v] = w-v+\varrho(w)-\varrho(v)$, we write $W_{2}$ as follows:
\[
W_{0} + V(g) - V(f) + O\biggl(\frac{x^{2}\ln{x}}{C}\biggr) + O\biggl(\frac{x^{2}(\ln{x})^{4}}{D}\biggr),
\]
where
\begin{multline*}
W_{0} = \sum\limits_{\Delta\leqslant D}\sum\limits_{\substack{0\leqslant\xi_{1},\xi_{2}< \Delta \\ \xi_{1}^{2}\equiv \xi_{2}^{2}\;(\mmod{\Delta})}}\prsum\limits_{\substack{1\leqslant m_{i}\leqslant\sqrt{x},\\ m_{i}\equiv \xi_{i}\;(\mmod{\Delta})\\ i = 1,2}}\prsum\limits_{\substack{\frac{\scriptstyle \sqrt{x}}{\scriptstyle C\Delta\mathstrut }<\kappa_{i}\leqslant \frac{\scriptstyle \sqrt{x}+m_{i}}{\scriptstyle \Delta\mathstrut },\; i=1,2 \\ (\kappa_{1},\kappa_{2})=1,\;\kappa_{1}\leqslant \kappa_{2}}}\frac{g-f}{\kappa_{2}},\\
V(\varphi) = \sum\limits_{\Delta\leqslant D}\sum\limits_{\substack{0\leqslant\xi_{1},\xi_{2}< \Delta \\ \xi_{1}^{2}\equiv \xi_{2}^{2}\;(\mmod{\Delta})}}\prsum\limits_{\substack{1\leqslant m_{i}\leqslant\sqrt{x},\\ m_{i}\equiv \xi_{i}\;(\mmod{\Delta}) \\ i = 1,2}}\prsum\limits_{\substack{\frac{\scriptstyle \sqrt{x}}{\scriptstyle C\Delta\mathstrut }<\kappa_{i}\leqslant \frac{\scriptstyle \sqrt{x}+m_{i}}{\scriptstyle \Delta\mathstrut },\; i=1,2 \\ (\kappa_{1},\kappa_{2})=1,\;\kappa_{1}\leqslant \kappa_{2}}}\varrho\biggl(\frac{\varphi - \nu\bar{\kappa}_{1}}{\kappa_{2}}\biggr)
\end{multline*}
(here $\varphi = f,g$). The sum $W_{0}$ contributes to the main term of $W$ and the sums $V(f), V(g)$ give the remainder term in the asymptotic formula for $W$.

\pagebreak

\textbf{\S 4. Asymptotic formula for $\boldsymbol{W_{0}}$.}
\vspace{0.5cm}

For brevity, here and below, in the sum signs, we write $\xi_{1}, \xi_{2}$ instead of $0\leqslant\xi_{1},\xi_{2}<\Delta$, $\xi_{1}^{2}\equiv \xi_{2}^{2}\;(\mmod{\Delta})$, and write $m_{1}, m_{2}$ instead of $1\leqslant m_{i}\leqslant\sqrt{x}$, $m_{i}\equiv \xi_{i}\;(\mmod{\Delta})$, $i = 1,2$. Thus we have
\begin{multline*}
W_{0} = \sum\limits_{\Delta\leqslant D}\;\sum\limits_{\xi_{1}, \xi_{2}}\;\prsum\limits_{m_{1}, m_{2}}\prsum\limits_{\substack{\frac{\scriptstyle \sqrt{x}}{\scriptstyle C\Delta\mathstrut } <\kappa_{i}\leqslant \frac{\scriptstyle \sqrt{x}+m_{i}}{\scriptstyle \Delta\mathstrut },\; i=1,2 \\ \kappa_{1}\leqslant \kappa_{2}}}\biggl(\,\sum\limits_{\delta | (\kappa_{1},\kappa_{2})}\mu(\delta)\biggr)\frac{g-f}{\kappa_{2}} =\\
= \sum\limits_{\Delta\leqslant D}\sum\limits_{1\leqslant \delta \leqslant \frac{\scriptstyle 2\sqrt{x}}{\scriptstyle \Delta}}\mu(\delta)\sum\limits_{\xi_{1}, \xi_{2}}
\prsum\limits_{m_{1}, m_{2}}
\prsum\limits_{\substack{\frac{\scriptstyle \sqrt{x}}{\scriptstyle C\Delta\mathstrut }<\kappa_{i}\leqslant \frac{\scriptstyle \sqrt{x}+m_{i}}{\scriptstyle \Delta\mathstrut }\\ \kappa_{i}\equiv 0\;(\mmod{\delta}),\; i=1,2 \\ \kappa_{1}\leqslant \kappa_{2}}}\frac{g-f}{\kappa_{2}}.
\end{multline*}
Next, we estimate the contribution coming to $W_{0}$ from $\delta>D$. Putting $\kappa_{i} = \delta n_{i}$ and using the inequalities
\begin{equation}\label{lab_12}
\frac{g-f}{\kappa_{2}}\leqslant \frac{g}{\kappa_{2}} < \frac{x}{q_{1}\kappa_{2}} = \frac{x}{\Delta\delta^{2}n_{1}n_{2}},
\end{equation}
we find that such contribution does not exceed
\begin{multline*}
x\sum\limits_{\Delta\leqslant D}\frac{1}{\Delta}\sum\limits_{D<\delta\leqslant \frac{\scriptstyle 2\sqrt{x}}{\scriptstyle \Delta}}\frac{1}{\delta^{2}}\sum\limits_{\xi_{1},\xi_{2}}
\sum\limits_{m_{1}, m_{2}}\sum\limits_{1\leqslant n_{i}\leqslant \frac{\scriptstyle 2\sqrt{x}}{\scriptstyle \Delta\delta}}\frac{1}{n_{1}n_{2}}\ll \\
\ll x(\ln{x})^{2}\sum\limits_{\Delta\leqslant D}\frac{1}{\Delta}\sum\limits_{D<\delta\leqslant \frac{\scriptstyle 2\sqrt{x}}{\scriptstyle \Delta}}\frac{\rho(\Delta)}{\delta^{2}}\biggl(\frac{\sqrt{x}}{\Delta}\biggr)^{2} \ll x^{2}(\ln{x})^{2}\sum\limits_{\Delta\leqslant D}\frac{\tau(\Delta)}{\Delta^{2}}\sum\limits_{\delta> D}\frac{1}{\delta^{2}}\ll \frac{x^{2}(\ln{x})^{2}}{D}.
\end{multline*}
If the tuple $(m_{1}, m_{2}, q_{1}, q_{2}) = (m_{1}, m_{2}, \Delta\delta n_{1}, \Delta\delta n_{2})$ is not in $\mathcal{A}(x)$ then at least one of the conditions $g>f$, $f>0$ fails. Hence, for such tuple we have $\max{(0,g-f)}= 0$. Otherwise, we obviously have $g-f = \max{(0,g-f)}$. Hence,
\begin{multline*}
W_{0} = \sum\limits_{\Delta\leqslant D}\sum\limits_{\delta\leqslant D}\frac{\mu(\delta)}{\delta}\sum\limits_{\xi_{1},\xi_{2}}W_{0}(\Delta,\delta;\xi_{1},\xi_{2}) + O\biggl(\frac{x^{2}(\ln{x})^{2}}{D}\biggr),\\
W_{0}(\Delta,\delta;\xi_{1},\xi_{2}) = \prsum\limits_{m_{1}, m_{2}}
\prsum\limits_{\substack{\frac{\scriptstyle \sqrt{x}}{\scriptstyle C\Delta\delta\mathstrut }<n_{i}\leqslant \frac{\scriptstyle \sqrt{x}+m_{i}}{\scriptstyle \Delta\delta\mathstrut },\; i=1,2 \\ n_{1}\leqslant n_{2}}}\frac{\max{(0,g-f)}}{n_{2}}.
\end{multline*}
Now we turn back to the representation of $f$ and $g$ in the form $f = \max{(f_{1},f_{2},f_{3},f_{4})}$, $g = \max{(g_{1}, g_{2}, g_{3})}$, and define $\mathcal{A}_{rs}(x)$ ($1\leqslant r\leqslant 4$, $1\leqslant s\leqslant 3$) as the subset of all quadruples $(m_{1}, m_{2}, q_{1}, q_{2}) = (m_{1}, m_{2}, \Delta\delta n_{1}, \Delta\delta n_{2})$ from $\mathcal{A}(x)$ such that
$f = f_{r}$, $g = g_{s}$. Let $\mathcal{A}_{rs}(x;\Delta,\delta;\xi_{1},\xi_{2}) = \mathcal{B}_{rs}$ be the set of all tuples from $\mathcal{A}_{rs}(x)$ that corresponds to given values $\Delta, \delta, \xi_{1}$ and $\xi_{2}$. Consequently, the sum $W_{0}(\Delta,\delta;\xi_{1},\xi_{2})$ splits into the sums $W_{rs} = W_{rs}(\Delta,\delta;\xi_{1},\xi_{2})$ corresponding to the tuples from the set $\mathcal{B}_{rs}$. Finally, we put $m_{i} = \Delta\ell_{i}+\xi_{i}$.

Note that any set $\mathcal{B}_{rs}$ is described by the system of inequalities which are linear with respect to variables $m_{1}, m_{2}, q_{1}, q_{2}$ (or $\ell_{1}, \ell_{2}, n_{1}, n_{2}$). Let us replace each sum $W_{rs}$ by the multiple integral. For brevity, we denote the variables of integration by the same letters: $\ell_{1}, \ell_{2}, n_{1}$ and  $n_{2}$.

Suppose that $W_{rs}$ has the form
\begin{multline*}
W_{rs} = \sum\limits_{(\ell_{1})}\sum\limits_{(\ell_{2})}\sum\limits_{(n_{1})}\sum\limits_{(n_{2})}\psi_{1}(\ell_{1}, \ell_{2}, n_{1}, n_{2}),\quad \text{where} \\
\psi_{1}(\ell_{1}, \ell_{2}, n_{1}, n_{2}) = \frac{1}{n_{2}}\bigl(g_{s}(\Delta\ell_{1}+\xi_{1},\Delta\ell_{2}+\xi_{2},\Delta\delta n_{1},\Delta\delta n_{2})-\\
-f_{r}(\Delta\ell_{1}+\xi_{1},\Delta\ell_{2}+\xi_{2},\Delta\delta n_{1},\Delta\delta n_{2})\bigr)
\end{multline*}
and $(\ell_{1}), (\ell_{2}), (n_{1}), (n_{2})$ are the systems of linear inequalities that define the domains of corresponding variables. Explicit formulas for $f_{r}$, $g_{s}$, $1\leqslant r\leqslant 4$, $1\leqslant s\leqslant 3$, imply that $\psi_{1}$ has the form
\begin{equation}\label{lab_13}
\frac{\mathscr{P}(\ell_{1},\ell_{2},n_{1},n_{2})}{n_{1}n_{2}},
\end{equation}
where $\mathscr{P}$ is a polynomial of degree at most two. Changing the order of summation, one may write $W_{rs}$ as follows:
\[
W_{rs} = \sum\limits_{(n_{2})^{*}}\sum\limits_{(n_{1})^{*}}\sum\limits_{(\ell_{2})^{*}}\sum\limits_{(\ell_{1})^{*}}\psi_{1}(\ell_{1}, \ell_{2}, n_{1}, n_{2}).
\]
Here $(n_{i})^{*},(\ell_{i})^{*}$ are linear systems of the form
\[
\mathscr{L}_{i}<\ell_{i}\leqslant \mathscr{L}_{i}^{'},\quad \mathscr{N}_{i}<n_{i}\leqslant \mathscr{N}_{i}^{'},\quad i=1,2,
\]
where $\mathscr{L}_{1},\mathscr{L}_{1}^{'}$ are linear forms in variables $\ell_{2}, n_{1}, n_{2}$;
$\mathscr{L}_{2},\mathscr{L}_{2}^{'}$ are linear forms in variables $n_{1}, n_{2}$;  $\mathscr{N}_{1},\mathscr{N}_{1}^{'}$ are linear forms in variable $n_{2}$, and
$\mathscr{N}_{2},\mathscr{N}_{2}^{'}$ depends only on $x$. At the same time, for all values of variables under considering we have
\[
0\leqslant \mathscr{L}_{i} < \mathscr{L}_{i}^{'} \leqslant \frac{\sqrt{x}}{\Delta},\quad \frac{\sqrt{x}}{C\Delta\delta}\leqslant \mathscr{N}_{i}<\mathscr{N}_{i}^{'}\leqslant \frac{\sqrt{x}}{\Delta\delta},\quad i = 1,2.
\]
In view of (\ref{lab_13}), for any fixed triple $n_{2}, n_{1},\ell_{2}$ the function $\psi_{1}$ (as a function of $\ell_{1}$) is a linear combination of piecewisely monotonic functions defined on $[\mathscr{L}_{1},\mathscr{L}_{1}^{'}]$. Obviously, both the number of terms in such linear combination and the number of segments of monotonicity of each term are bounded above by some absolute constant. The application of Lemma 1 together with (12) yields:
\begin{multline}\label{lab_14}
\sum\limits_{(\ell_{1})^{*}}\psi_{1}(\ell_{1}, \ell_{2}, n_{1}, n_{2}) = \psi_{2}(\ell_{2}, n_{1}, n_{2}) + O\biggl(\frac{x}{\Delta\delta^{2} n_{1}n_{2}}\biggr),\\
\psi_{2}(\ell_{2}, n_{1}, n_{2}) = \int_{(\ell_{1})^{*}}\psi_{1}(\ell_{1}, \ell_{2}, n_{1}, n_{2})d\ell_{1}.
\end{multline}
The contribution coming to $W_{rs}$ from the remainder term of (\ref{lab_14}) does not exceed
\[
\frac{x}{\Delta\delta^{2}}\sum\limits_{\substack{1\leqslant n_{i}\leqslant \sqrt{x}/(\Delta\delta) \\ i = 1,2}}\frac{1}{n_{1}n_{2}}\sum\limits_{0\leqslant \ell_{2}\leqslant \sqrt{x}/\Delta}1\ll \frac{x^{3/2}(\ln{x})^{2}}{\Delta^{2}\delta^{2}}.
\]
Thus we have
\[
W_{rs} =  \sum\limits_{(n_{2})^{*}}\sum\limits_{(n_{1})^{*}}\sum\limits_{(\ell_{2})^{*}}\psi_{2}(\ell_{2}, n_{1}, n_{2}) + O\biggl(\frac{x^{3/2}(\ln{x})^{2}}{\Delta^{2}\delta^{2}}\biggr).
\]
Further, it follows from (\ref{lab_13}) that $\psi_{2}(\ell_{2}, n_{1}, n_{2})$ has the form
\begin{equation}\label{lab_15}
\frac{\mathscr{Q}(\ell_{2},n_{1},n_{2})}{n_{1}n_{2}},
\end{equation}
where $\mathscr{Q}$ is a polynomial of degree at most three. Moreover, the inequality (\ref{lab_12}) implies that
\begin{equation}\label{lab_16}
\psi_{2}(\ell_{2}, n_{1}, n_{2}) \ll \int_{(\ell_{1})^{*}}\frac{x\,d\ell_{1}}{\Delta\delta^{2}n_{1}n_{2}}\ll \frac{x^{3/2}}{\Delta^{2}\delta^{2}n_{1}n_{2}}.
\end{equation}
In view of (\ref{lab_15}), for any $n_{1}, n_{2}$, the function $\psi_{2}$ (as a function of $\ell_{2}$) is a linear combination of piecewisely monotonic functions defined on $[\mathscr{L}_{2},\mathscr{L}_{2}^{'}]$. Both the number of terms in such linear combination and the number of segments of monotonicity of each term are bounded above by some absolute constant. Using Lemma 1 and (\ref{lab_16}) we get
\begin{multline}\label{lab_17}
\sum\limits_{(\ell_{2})^{*}}\psi_{2}(\ell_{2}, n_{1}, n_{2}) = \psi_{3}(n_{1}, n_{2}) + O\biggl(\frac{x^{3/2}}{\Delta^{2}\delta^{2} n_{1}n_{2}}\biggr),\\
\psi_{3}(n_{1}, n_{2}) = \int_{(\ell_{2})^{*}}\psi_{2}(\ell_{2}, n_{1}, n_{2})d\ell_{2} = \int_{(\ell_{2})^{*}}\int_{(\ell_{1})^{*}}\psi_{1}(\ell_{1},\ell_{2}, n_{1}, n_{2})d\ell_{1}d\ell_{2}.
\end{multline}
The contribution coming to $W_{rs}$ from $O$-term of (\ref{lab_17}) is less than
\[
\frac{x^{3/2}}{\Delta^{2}\delta^{2}}\sum\limits_{1\leqslant n_{1}, n_{2}\leqslant \sqrt{x}/(\Delta\delta)}\frac{1}{n_{1}n_{2}}\ll \frac{x^{3/2}(\ln{x})^{2}}{\Delta^{2}\delta^{2}}.
\]
Hence,
\[
W_{rs} =  \sum\limits_{(n_{2})^{*}}\sum\limits_{(n_{1})^{*}}\psi_{3}(n_{1}, n_{2}) + O\biggl(\frac{x^{3/2}(\ln{x})^{2}}{\Delta^{2}\delta^{2}}\biggr).
\]
From (\ref{lab_15}), it follows that $\psi_{3}(n_{1}, n_{2})$ has the form
\begin{equation}\label{lab_18}
\frac{\mathscr{R}(n_{1},n_{2})}{n_{1}n_{2}},
\end{equation}
where $\mathscr{R}$ is a polynomial of degree at most four. The inequality (\ref{lab_16}) implies that
\begin{equation}\label{lab_19}
\psi_{3}(n_{1}, n_{2}) \ll \int_{(\ell_{2})^{*}}\psi_{2}(\ell_{2}, n_{1}, n_{2})d\ell_{2} \ll \frac{x^{2}}{\Delta^{3}\delta^{2}n_{1}n_{2}}.
\end{equation}
Since $\sqrt{x}/(C\Delta\delta)\leqslant n_{1}\leqslant \sqrt{x}/(\Delta\delta)$ then
\begin{equation}\label{lab_20}
\psi_{3}(n_{1}, n_{2}) \ll \frac{Cx^{3/2}}{\Delta^{2}\delta n_{2}}
\end{equation}
on the segment $[\mathscr{N}_{1},\mathscr{N}_{1}^{'}]$. Applying Lemma 1 to the sum over $n_{1}$ we find that
\begin{multline}\label{lab_21}
\sum\limits_{(n_{1})^{*}}\psi_{3}(n_{1}, n_{2}) = \psi_{4}(n_{2}) + O\biggl(\frac{Cx^{3/2}}{\Delta^{2}\delta n_{2}}\biggr),\quad \text{where} \\
\psi_{4}(n_{2}) = \int_{(n_{1})^{*}}\psi_{3}(n_{1}, n_{2})dn_{1} = \int_{(n_{1})^{*}}\int_{(\ell_{2})^{*}}\int_{(\ell_{1})^{*}}\psi_{1}(\ell_{1},\ell_{2}, n_{1}, n_{2})d\ell_{1}d\ell_{2}dn_{1}.
\end{multline}
The contribution coming to $W_{rs}$ from $O$-term of (\ref{lab_21}) does not exceed in order
\[
\frac{Cx^{3/2}}{\Delta^{2}\delta}\sum\limits_{1\leqslant n_{2}\leqslant \sqrt{x}/(\Delta\delta)}\frac{1}{n_{2}}\ll \frac{Cx^{3/2}\ln{x}}{\Delta^{2}\delta}.
\]
Therefore,
\[
W_{rs} = \sum\limits_{(n_{2})^{*}}\psi_{4}(n_{2}) + O\biggl(\frac{x^{3/2}(\ln{x})^{2}}{\Delta^{2}\delta^{2}}\biggr) + O\biggl(\frac{Cx^{3/2}\ln{x}}{\Delta^{2}\delta}\biggr)
= \sum\limits_{(n_{2})^{*}}\psi_{4}(n_{2}) + O\biggl(\frac{Cx^{3/2}\ln{x}}{\Delta^{2}\delta}\biggr).
\]
Moreover, the estimate (\ref{lab_20}) implies that
\begin{equation*}
\psi_{4}(n_{2}) \ll \int_{(n_{1})^{*}}\frac{x^{2}}{\Delta^{3}\delta^{2}}\,\frac{dn_{1}}{n_{1}n_{2}} \ll \frac{x^{2}\ln{x}}{\Delta^{3}\delta^{2}n_{2}} \ll \frac{Cx^{3/2}\ln{x}}{\Delta^{2}\delta}.
\end{equation*}
Finally, in view of (\ref{lab_18}), $\psi_{4}(n_{2})$ has the form
\begin{equation}\label{lab_22}
\frac{1}{n_{2}}\bigl(\mathscr{S}_{1}(n_{1}) + \mathscr{S}_{2}(n_{2})\ln{\mathscr{N}}\bigr),
\end{equation}
where $\mathscr{S}_{j}(n_{2})$ are polynomials of degree at most three and $\mathscr{N}$ is one of the forms $\mathscr{N}_{1}$, $\mathscr{N}_{1}^{'}$. Thus, $\psi_{4}$ is a linear combination of
of piecewisely monotonic functions defined on $[\mathscr{N}_{2},\mathscr{N}_{2}^{'}]$. Both the number of terms in such linear combination and the number of segments of monotonicity of each term are bounded above by some absolute constant. By Lemma 1 and (22), we find that
\begin{multline*}
\sum\limits_{(n_{2})^{*}}\psi_{4}(n_{2}) = \int_{(n_{2})^{*}}\psi_{4}(n_{2})dn_{2} + O\biggl(\frac{Cx^{3/2}\ln{x}}{\Delta^{2}\delta}\biggr), \\
W_{rs} = \int_{(n_{2})^{*}}\psi_{4}(n_{2})dn_{2} + O\biggl(\frac{Cx^{3/2}\ln{x}}{\Delta^{2}\delta}\biggr) = \\
= \int_{(n_{2})^{*}}\int_{(n_{1})^{*}}\int_{(\ell_{2})^{*}}\int_{(\ell_{1})^{*}}\psi_{1}(\ell_{1},\ell_{2}, n_{1}, n_{2})d\ell_{1}d\ell_{2}dn_{1}dn_{2} + O\biggl(\frac{Cx^{3/2}\ln{x}}{\Delta^{2}\delta}\biggr) = \\
= \int_{(\ell_{1})}\int_{(\ell_{2})}\int_{(n_{1})}\int_{(n_{2})}\psi_{1}(\ell_{1},\ell_{2}, n_{1}, n_{2})d\ell_{1}d\ell_{2}dn_{1}dn_{2} + O\biggl(\frac{Cx^{3/2}\ln{x}}{\Delta^{2}\delta}\biggr).
\end{multline*}
The last integral is expressed as follows:
\[
\int_{(\ell_{1})}\int_{(\ell_{2})}\int_{(n_{1})}\int_{(n_{2})}\chi_{rs}\frac{g_{s}-f_{r}}{n_{2}}\,dn_{2}dn_{1}d\ell_{2}d\ell_{1} =  \iiiint\limits_{\Omega}\frac{\chi_{rs}}{n_{2}}\,\max{(0,g-f)}\,dn_{2}dn_{1}d\ell_{2}d\ell_{1},
\]
where $\chi_{rs} = \chi_{rs}(\ell_{1},\ell_{2},n_{1},n_{2})$ is the characteristic function of the set of points $(\ell_{1},\ell_{2},n_{1},n_{2})\in \mathbb{R}_{+}^{4}$
satisfying the conditions $f = f_{r}$, $g = g_{s}$, $f_{r}<g_{s}$, and $\Omega = \Omega(\Delta,\delta;\xi_{1},\xi_{2})$ denotes the set of points $(\ell_{1},\ell_{2},n_{1},n_{2})\in \mathbb{R}_{+}^{4}$ obeying the following set of conditions:
\begin{equation*}
\begin{cases}
0<m_{i}\leqslant\sqrt{x},\quad m_{i} = \Delta\ell_{i}+\xi_{i},\quad i = 1,2, \\
\frac{\displaystyle \sqrt{x}}{\displaystyle C\Delta\delta\mathstrut} < n_{i}\leqslant \frac{\displaystyle \sqrt{x}+m_{i}}{\displaystyle \Delta\delta\mathstrut},\quad i = 1,2, \\
n_{1}\leqslant n_{2}.
\end{cases}
\end{equation*}
The summation over all pairs $r,s$ yields:
\[
W_{0}(\Delta,\delta;\xi_{1},\xi_{2}) = I_{0} + O\biggl(\frac{Cx^{3/2}\ln{x}}{\Delta^{2}\delta}\biggr),\quad I_{0} = \iiiint\limits_{\Omega}\frac{\max{(0,g-f)}}{n_{2}}\,dn_{2}dn_{1}d\ell_{2}d\ell_{1}.
\]
The change of variables $m_{i} = \Delta\ell_{i}+\xi_{i}, q_{i} = \Delta\delta n_{i}$ reduces this integral to
\begin{multline*}
\frac{\Delta\delta}{\Delta^{4}\delta^{2}}\int_{0}^{\sqrt{x}}dm_{1}\int_{0}^{\sqrt{x}}dm_{2}\iint\limits_{\substack{\sqrt{x}/C<q_{i}\leqslant\sqrt{x}+m_{i},\,i=1,2 \\ q_{1}\leqslant q_{2}}}\psi(m_{1},m_{2},q_{1},q_{2})dq_{1}dq_{2},\\
\psi(m_{1},m_{2},q_{1},q_{2}) = \\
=\frac{1}{q_{1}q_{2}}\,\max{\bigl(0,\min{\{x,(q_{1}+m_{1})^{2},(q_{2}+m_{2})^{2}\}}-\max{\{m_{1}^{2},m_{2}^{2},(q_{1}-m_{1})^{2},(q_{2}-m_{2})^{2}\}}\bigr)}.
\end{multline*}
Putting $m_{i} = u_{i}\sqrt{x}$, $q_{i} = v_{i}\sqrt{x}$, $i = 1,2$, we get
\begin{multline*}
I_{0} = \frac{x^{2}}{2\Delta^{3}\delta}\,\mathfrak{S}\biggl(\frac{1}{C}\biggr),\quad \mathfrak{S}(\varepsilon) =
2\int_{0}^{1}du_{1}\int_{0}^{1}du_{2}\iint\limits_{\substack{\varepsilon\leqslant v_{i}\leqslant 1+u_{i},\;i=1,2 \\ v_{1}\leqslant v_{2}}}\varphi(u_{1},u_{2},v_{1},v_{2})dv_{1}dv_{2} = \\
= \int_{0}^{1}du_{1}\int_{0}^{1}du_{2}\int_{\varepsilon}^{1+u_{1}}dv_{1}\int_{\varepsilon}^{1+u_{2}}\varphi(u_{1},u_{2},v_{1},v_{2})dv_{2}, \\
\varphi(u_{1},u_{2},v_{1},v_{2}) = \\
=  \frac{1}{v_{1}v_{2}}\max{\bigl\{0,\min{(1,(u_{1}+v_{1})^{2},(u_{2}+v_{2})^{2})} - \max{(u_{1}^{2},u_{2}^{2},(u_{1}-v_{1})^{2},(u_{2}-v_{2})^{2})}\bigr\}}.
\end{multline*}
One can check (see \S 6) that
\[
\mathfrak{S}(\varepsilon) = \mathfrak{S}(0) + O\biggl(\varepsilon\ln^{2}\frac{1}{\varepsilon}\biggr).
\]
Setting $\mathfrak{S} = \mathfrak{S}(0)$, we find
\[
I_{0} = \frac{x^{2}}{\Delta^{3}\delta}\biggl(\mathfrak{S} + O\biggl(\frac{(\ln{C})^{2}}{C}\biggr)\biggr),\quad W_{0}(\Delta,\delta;\xi_{1},\xi_{2}) = \frac{x^{2}\mathfrak{S}}{\Delta^{3}\delta}
+ O\biggl(\frac{x^{2}(\ln{x})^{2}}{C\Delta^{3}\delta}\biggr)+ O\biggl(\frac{Cx^{3/2}\ln{x}}{\Delta^{2}\delta}\biggr).
\]
The summation over $\xi_{1},\xi_{2},\delta$ and $\Delta$ yields:
\begin{multline*}
W_{0} = 2\sum\limits_{\Delta\leqslant D}\sum\limits_{\delta\leqslant D}\frac{\mu(\delta)}{\delta}\rho(\Delta)\biggl\{\frac{x^{2}\mathfrak{S}}{2\Delta^{3}\delta}
+ O\biggl(\frac{x^{2}(\ln{x})^{2}}{C\Delta^{3}\delta}\biggr)+ O\biggl(\frac{Cx^{3/2}\ln{x}}{\Delta^{2}\delta}\biggr)\biggr\} + O\biggl(\frac{x^{2}(\ln{x})^{2}}{D}\biggr) = \\
= \frac{6}{\pi^{2}}\biggl(\sum\limits_{\Delta=1}^{+\infty}\frac{\rho(\Delta)}{\Delta^{3}}\biggr)\mathfrak{S}x^{2} + O\biggl(\frac{x^{2}(\ln{x})^{2}}{D}\biggr)
+ O\biggl(\frac{x^{2}(\ln{x})^{2}}{C}\biggr) + O\bigl(Cx^{3/2}(\ln{x})^{2}\bigr).
\end{multline*}
Taking $C = D = \sqrt[4\;]{x}/2$ and using Lemma 3, we obtain
\[
W_{0} = c_{0}x^{2} + O\bigl(x^{7/4}(\ln{x})^{2}\bigr),\quad c_{0} = \frac{13}{14}\,\frac{\zeta(2)}{\zeta(3)}\,\mathfrak{S}.
\]

\textbf{\S 5. Estimates of the sums $\boldsymbol{V(f), V(g)}$.}
\vspace{0.5cm}

Obviously, the contribution coming to $V(\varphi)$ from the pairs $m_{1} = m_{2}$ is less than
\[
\sum\limits_{1\leqslant \Delta \leqslant D}\sum\limits_{\xi_{1}=1}^{\Delta}\sum\limits_{\substack{1\leqslant m_{1}\leqslant \sqrt{x} \\ m_{1}\equiv \xi_{1}\;(\mmod{\Delta})}}\biggl(\frac{\sqrt{x}}{\Delta}\biggr)^{2}\ll \sum\limits_{1\leqslant \Delta \leqslant D}\Delta\biggl(\frac{\sqrt{x}}{\Delta}\biggr)^{3}\ll x^{3/2}.
\]
So, in what follows, we assume that $m_{1}\ne m_{2}$.

Thus we write $V(\varphi)$ in the form
\[
\sum\limits_{1\leqslant \Delta \leqslant D}\;\sum\limits_{\substack{0<\xi_{1},\xi_{2}\leqslant \Delta \\ \xi_{1}^{2}\equiv \xi_{2}^{2}\;(\mmod{\Delta})}}V(\varphi,\Delta;\xi_{1},\xi_{2}) + O(x^{3/2}),
\]
where
\begin{multline*}
V(\varphi,\Delta;\xi_{1},\xi_{2}) = \prsum\limits_{\substack{1\leqslant m_{i}\leqslant \sqrt{x} \\ m_{i}\equiv\xi_{i}\;(\mmod{\Delta}) \\ i = 1,2}}V(\varphi,\Delta;\xi_{1},\xi_{2};m_{1},m_{2}),\\
V(\varphi,\Delta;\xi_{1},\xi_{2};m_{1},m_{2}) = \prsum\limits_{\substack{\frac{\scriptstyle\sqrt{x}}{\scriptstyle C\Delta\mathstrut}<\kappa_{i}\leqslant \frac{\scriptstyle\sqrt{x}}{\scriptstyle\Delta\mathstrut},\;  i=1,2 \\ (\kappa_{1},\kappa_{2})=1,\;\kappa_{1}\leqslant \kappa_{2}}}\varrho\biggl(\frac{\varphi - \nu\overline{\kappa}_{1}}{\kappa_{2}}\biggr).
\end{multline*}
Suppose now that $H>1$ is integer whose precise value will be chosen later, and denote by $U^{(1)}$ and $U^{(2)}$ the contributions coming to the last sum from the terms $\varrho_{H}$ and $r_{H}$ of Lemma 2. Arguing as above, we split $U^{(j)}$ to the sums $U_{rs}^{(j)}$, $j=1,2$, where the indexes $r,s$ mean that the corresponding tuples $(m_{1},m_{2},q_{1},q_{2}) = (m_{1}, m_{2},\Delta\kappa_{1},\Delta\kappa_{2})$ belong to the set $\mathcal{A}_{rs}(x)$. Then, for fixed $r,s$ we write
\[
U_{rs}^{(1)} = \sum\limits_{(\kappa_{1})}\sum\limits_{(\kappa_{2})}\varrho_{H}\biggl(\frac{\varphi - \nu\overline{\kappa}_{1}}{\kappa_{2}}\biggr),
\]
where the symbols $(\kappa_{1})$ and $(\kappa_{2})$ mean the domains of the type $\sigma_{1}<\kappa_{1}\leqslant \tau_{1}$, $\sigma_{2}<\kappa_{2}\leqslant \tau_{2}$, $(\kappa_{1},\kappa_{1})=1$ and
\[
\sigma_{1} = \sigma_{1}(x;m_{1},m_{2}),\quad \tau_{1} = \tau_{1}(x;m_{1},m_{2}),\quad \sigma_{2} = \sigma_{2}(x;m_{1},m_{2},\kappa_{1}),\quad \tau_{2} = \tau_{2}(x;m_{1},m_{2},\kappa_{1})
\]
are linear functions of $m_{1},m_{2},\kappa_{1}$. The restrictions imposed on $\kappa_{i}$ imply the inequalities $E<\kappa_{i}\leqslant K$, where
\[
E = \biggl[\frac{\sqrt{x}}{2C\Delta\mathstrut}\biggr],\quad K = 2\biggl[\frac{\sqrt{x}}{\Delta\mathstrut}\biggr]+3.
\]
Further, using Lemma 2 we have
\begin{multline*}
U_{rs}^{(1)} = \frac{1}{2\pi i}\sum\limits_{0<|h|\leqslant H}\frac{1}{h}\sum\limits_{(\kappa_{1})}\sum\limits_{(\kappa_{2})}e\biggl(\frac{h\varphi - h\nu\bar{\kappa}_{1}}{\kappa_{2}}\biggr)
= \frac{1}{2\pi i}\prsum\limits_{L\leqslant H/2}\frac{U_{rs}^{(1)}(L)}{L},\\
U_{rs}^{(1)}(L) = \sum\limits_{L<|h|\leqslant 2L}\gamma(h)\sum\limits_{(\kappa_{1})}\sum\limits_{(\kappa_{2})}e\biggl(\frac{h\varphi - h\nu\bar{\kappa}_{1}}{\kappa_{2}}\biggr),\quad \gamma(h) = \frac{L}{h}.
\end{multline*}
Next, we transform the domain  $(\kappa_{2})$ to make the bounds for $\kappa_{2}$ independent of $\kappa_{1}$. Thus we get
\begin{multline*}
\sum\limits_{(\kappa_{2})}e\biggl(\frac{h\varphi - h\nu\bar{\kappa}_{1}}{\kappa_{2}}\biggr) = \sum\limits_{\substack{\sigma_{2}<\kappa_{2}\leqslant \tau_{2} \\ (\kappa_{1},\kappa_{2})=1}}
e\biggl(\frac{h\varphi - h\nu\bar{\kappa}_{1}}{\kappa_{2}}\biggr) = \\
=\, \sum\limits_{\substack{E<\kappa_{2}\leqslant K \\ (\kappa_{1},\kappa_{2})=1}}\frac{1}{K}\sum\limits_{|a|<K/2}\sum\limits_{\sigma_{2}<\lambda\leqslant \tau_{2}}e\biggl(\frac{a}{K}(\kappa_{2}-\lambda)\biggr)e\biggl(\frac{h\varphi - h\nu\bar{\kappa}_{1}}{\kappa_{2}}\biggr) =\\
=\,\sum\limits_{|a|<K/2}\frac{1}{|a|+1}\sum\limits_{\substack{E<\kappa_{2}\leqslant K \\ (\kappa_{1},\kappa_{2})=1}}\alpha(\kappa_{2})\beta(\kappa_{1})
e\biggl(\frac{h\varphi - h\nu\bar{\kappa}_{1}}{\kappa_{2}}\biggr),
\end{multline*}
where
\[
\alpha(\kappa_{2}) = e\biggl(\frac{a\kappa_{2}}{K}\biggr),\quad \beta(\kappa_{1}) = \beta_{a}(\kappa_{1}) = \frac{|a|+1}{K}\sum\limits_{\sigma_{2}<\lambda\leqslant \tau_{2}}e\biggl(-\frac{a\lambda}{K}\biggr).
\]
It is easy to check that $|\beta_{a}(\kappa_{1})|\leqslant 1$ for any $a$. Hence,
\begin{multline*}
U_{rs}^{(1)}(L) = \sum\limits_{|a|<K/2}\frac{U_{rs}^{(1)}(L;a)}{|a|+1},\quad\text{where}\\
U_{rs}^{(1)}(L;a) = \sum\limits_{L<|h|\leqslant 2L}\gamma(h)\sum\limits_{(\kappa_{1})}\beta(\kappa_{1})\sum\limits_{\substack{E<\kappa_{2}\leqslant K \\ (\kappa_{1},\kappa_{2})=1}}\alpha(\kappa_{2})e\biggl(\frac{h\varphi - h\nu\bar{\kappa}_{1}}{\kappa_{2}}\biggr) =\\
=\,\sum\limits_{L<|h|\leqslant 2L}\sum\limits_{\sigma_{1}<\kappa_{1}\leqslant \tau_{1}}\sum\limits_{\substack{E<\kappa_{2}\leqslant K \\ (\kappa_{1},\kappa_{2})=1}}\alpha(\kappa_{2})\beta(\kappa_{1})\gamma(h)e\biggl(\frac{h\varphi - h\nu\bar{\kappa}_{1}}{\kappa_{2}}\biggr).
\end{multline*}
Finally, we split the domains of $\kappa_{1}$  and $\kappa_{2}$ into the segments of the type $M<\kappa_{1}\leqslant M_{1}$, $N<\kappa_{2}\leqslant N_{1}$, where $M_{1} = 2M$ for all cases except at most one, where $M_{1}<2M$ (similarly for $N$, $N_{1}$).

Limiting ourselves by the case of positive $h$ and taking the conjugation we lead to the sums
\[
U(L,M,N) = \sum\limits_{L<h\leqslant 2L}\mathop{\sum_{M<\kappa_{1}\leqslant M_{1}}\sum_{N<\kappa_{2}\leqslant N_{1}}}\limits_{(\kappa_{1},\kappa_{2})=1}
\alpha(\kappa_{2})\beta(\kappa_{1})\gamma(h)e\biggl(\frac{\nu h\overline{\kappa}_{1}}{\kappa_{2}} + \mathscr{F}(\kappa_{1},\kappa_{2})\biggr),
\]
where $\mathscr{F}(\kappa_{1},\kappa_{2}) = -h\varphi/\kappa_{2}$. Now we have to estimate the derivatives $\partial \mathscr{F}/\partial\kappa_{i}$, $i=1,2$.
Thus, in the case
\[
\varphi = f_{3} = \frac{(q_{2}-m_{2})^{2}-m_{1}^{2}}{q_{1}} = \frac{(\Delta\kappa_{2}-m_{2})^{2} - m_{1}^{2}}{\Delta\kappa_{1}}
\]
we find
\[
\mathscr{F} = \frac{h(m_{1}^{2}-(\Delta\kappa_{2}-m_{2})^{2})}{\Delta\kappa_{1}\kappa_{2}},\quad \frac{\partial \mathscr{F}}{\partial\kappa_{1}} = \frac{h((\Delta\kappa_{2}-m_{2})^{2}-m_{1}^{2})}{\Delta\kappa_{1}^{2\mathstrut}\kappa_{2}},\quad \frac{\partial \mathscr{F}}{\partial\kappa_{2}} =
\frac{h(m_{2}^{2}-m_{1}^{2}-(\Delta\kappa_{2})^{2})}{\Delta\kappa_{1}\kappa_{2}^{2}}
\]
and therefore
\[
\biggl|\frac{\partial \mathscr{F}}{\partial\kappa_{1}}\biggr|\leqslant \frac{h(\sqrt{x})^{2}}{\Delta\kappa_{1}^{2}\kappa_{2}}\leqslant \frac{2xL}{\Delta\kappa_{1}^{2}\kappa_{2}},\quad
\biggl|\frac{\partial \mathscr{F}}{\partial\kappa_{2}}\biggr|\leqslant \frac{h\bigl((\sqrt{x})^{2}+(2\sqrt{x})^{2}\bigr)}{\Delta\kappa_{1}^{2}\kappa_{2}} \leqslant \frac{10xL}{\Delta\kappa_{1}^{2}\kappa_{2}}.
\]
Thus, the inequalities (\ref{lab_04}) hold with $X\ll xL/\Delta$. The same estimates are valid in all other cases. Setting $\vartheta = \nu = (m_{2}^{2}-m_{1}^{2})/\Delta$ in Lemma 4, we conclude that
\[
\frac{|\vartheta|L+X}{MN}\ll \frac{xL}{\Delta MN},\quad T = \biggl(1+\frac{|\vartheta|L+X}{MN}\biggr)^{1/2}\ll 1 + \sqrt{\frac{xL}{\Delta MN}}.
\]
Next, taking $\{\alpha(\kappa_{1})\}$, $\{\beta(\kappa_{2})\}$ and $\{\gamma(h)\}$ as the sequences $\mathbf{a}=\{a_{m}\}$, $\mathbf{b}=\{b_{n}\}$ and $\mathbf{c}=\{c_\ell\}$ in Lemma 4, we get $\|\mathbf{a}\|\ll\sqrt{M}$,
$\|\mathbf{b}\|\ll\sqrt{N}$, $\|\mathbf{c}\|\ll\sqrt{L}$. Hence, for any fixed $\varepsilon>0$, we obtain:
\begin{multline*}
U(L,M,N)\ll (LMN)^{\varepsilon/4}(LMN)^{1/2}\biggl(1+ \sqrt{\frac{xL}{\Delta MN}}\biggr)\times\\
\times\bigl\{(LMN)^{7/20}(M+N)^{1/4}+(LMN)^{3/8}(M+N)^{1/8}L^{1/8}\bigr\} \ll \\
\ll\,x^{3\varepsilon/4}\biggl(1+ \sqrt{\frac{xL}{\Delta MN}}\biggr)\bigl\{(LMN)^{17/20}(M+N)^{1/4}+(LMN)^{7/8}(M+N)^{1/8}L^{1/8}\bigr\}.
\end{multline*}
Now we set $M = K\cdot 2^{-\alpha}$, $N = K\cdot 2^{-\beta}$, where $1\leqslant 2^{\alpha}, 2^{\beta}\ll C$. Then we get
\begin{multline*}
\frac{xL}{\Delta MN}\ll \frac{xL}{\Delta}\cdot\frac{2^{\alpha+\beta}}{K^{2\mathstrut}} \ll \frac{xL}{\Delta}\cdot\frac{\Delta^{2}}{x}\,2^{\alpha+\beta}\ll \Delta L\cdot 2^{\alpha+\beta},\\
1 + \sqrt{\frac{xL}{\Delta MN}}\ll \sqrt{\Delta L}\cdot 2^{(\alpha + \beta)/2},\quad LMN \ll \frac{xL}{\Delta^{2}}\cdot 2^{-(\alpha+\beta)},\\
(LMN)^{17/20}(M+N)^{1/4} \ll \biggl(\frac{x}{\Delta^{2}}\biggr)^{39/40}L^{17/20}\bigl( 2^{-11\alpha/10-17\beta/20} + 2^{-17\alpha/20-11\beta/10}\bigr),\\
(LMN)^{7/8}(M+N)^{1/8}L^{1/8} \ll \biggl(\frac{x}{\Delta^{2}}\biggr)^{15/16}L\bigl(2^{-\alpha-7\beta/8} + 2^{-7\alpha/8-\beta}\bigr)
\end{multline*}
and hence
\begin{multline*}
U(L,M,N)\ll x^{3\varepsilon/4}\sqrt{\Delta L}\biggl\{\biggl(\frac{x}{\Delta^{2}}\biggr)^{39/40}L^{17/20}\bigl( 2^{-3\alpha/5-7\beta/20} + 2^{-7\alpha/20-3\beta/5}\bigr) + \\
+ \biggl(\frac{x}{\Delta^{2}}\biggr)^{15/16}L\bigl(2^{-\alpha/2-3\beta/8} + 2^{-3\alpha/8-\beta/2}\bigr)\biggr\}.
\end{multline*}
The summation over all $\alpha, \beta$ yields:
\begin{multline*}
U_{rs}^{(1)}(L,a) \ll x^{3\varepsilon/4}\biggl\{ \frac{x^{39/40}}{\Delta^{29/20}}\,L^{27/20} + \frac{x^{15/16}}{\Delta^{11/8}}\,L^{3/2} \biggr\},\\
U_{rs}^{(1)}(L) \ll x^{4\varepsilon/5}\biggl\{ \frac{x^{39/40}}{\Delta^{29/20}}\,L^{27/20} + \frac{x^{15/16}}{\Delta^{11/8}}\,L^{3/2} \biggr\}.
\end{multline*}
Further, taking $L = H\cdot 2^{-\gamma}$ we obtain
\[
U_{rs}^{(1)} \ll \sum_{\gamma}\frac{|U_{rs}^{(1)}(L)|}{L}\ll x^{5\varepsilon/6}\biggl\{\frac{x^{39/40}}{\Delta^{29/20}}\,H^{7/20} + \frac{x^{15/16}}{\Delta^{11/8}}\,H^{1/2} \biggr\}.
\]
Finally, summing over all pairs $r,s$ we find that
\begin{equation}\label{lab_23}
U^{(1)} \ll  x^{5\varepsilon/6}\biggl\{\frac{x^{39/40}}{\Delta^{29/20}}\,H^{7/20} + \frac{x^{15/16}}{\Delta^{11/8}}\,H^{1/2} \biggr\}.
\end{equation}
Now we estimate the sum $U^{(2)}$. By Lemma 2, for any $r,s$ and for some absolute constant $c>0$ we have
\begin{multline*}
|U_{rs}^{(2)}| = \biggl|\sum\limits_{(\kappa_{1})}\sum\limits_{(\kappa_{2})}r_{H}\biggl(\frac{\varphi - \nu\overline{\kappa}_{1}}{\kappa_{2}}\biggr)\biggr|\leqslant
c\sum\limits_{(\kappa_{1})}\sum\limits_{(\kappa_{2})}\vartheta_{H}\biggl(\frac{\varphi - \nu\overline{\kappa}_{1}}{\kappa_{2}}\biggr)\leqslant \\
\leqslant c\sum\limits_{\substack{E<\kappa_{1},\kappa_{2}\leqslant K \\ (\kappa_{1},\kappa_{2})=1}}\vartheta_{H}\biggl(\frac{\varphi - \nu\overline{\kappa}_{1}}{\kappa_{2}}\biggr) = c\sum\limits_{h = -\infty}^{+\infty}C(h)\sum\limits_{\substack{E<\kappa_{1},\kappa_{2}\leqslant K \\ (\kappa_{1},\kappa_{2})=1}}e\biggl(\frac{h\varphi - h\nu\overline{\kappa}_{1}}{\kappa_{2}}\biggr).
\end{multline*}
Let $\varepsilon_{1} = \varepsilon/10$ and $P = \bigl[H^{1+\varepsilon_{1}}\bigr]$. Taking $A=3$ in Lemma 2, we estimate the tail with $|h|>P$ as follows:
\[
c\sum\limits_{|h|>P}|C(h)|K^{2}\ll K^{2}H^{-3}\ll  \frac{x}{\Delta^{2}H^{3}}.
\]
Further, the term with $h = 0$ contributes at most
\[
c|C(0)|K^{2}\ll \frac{K^{2}\ln{H}}{H}\ll \frac{x\ln{x}}{\Delta^{2}H}.
\]
Now it remains to consider the contribution coming from $0<|h|\leqslant P$. To do this, we split the segment $1\leqslant h\leqslant P$ into the segments of the type $L<h\leqslant L_{1}$, $L_{1}\leqslant 2L$. Thus we find:
\begin{multline*}
U_{rs}(L) = \sum\limits_{L<h\leqslant L_{1}}C(h)\sum\limits_{\substack{E<\kappa_{1},\kappa_{2}\leqslant K \\ (\kappa_{1},\kappa_{2})=1}}e\biggl(\frac{h\varphi - h\nu\overline{\kappa}_{1}}{\kappa_{2}}\biggr) =\\
=\,\frac{\ln{H}}{H}\sum\limits_{L<h\leqslant L_{1}}
\mathop{\sum\limits_{E<\kappa_{1}\leqslant K}\sum\limits_{E<\kappa_{2}\leqslant K}}\limits_{(\kappa_{1},\kappa_{2})=1}
\gamma(h)e\biggl(\frac{h\varphi - h\nu\overline{\kappa}_{1}}{\kappa_{2}}\biggr),
\end{multline*}
where $\gamma(h) = HC(h)/\ln{H}\ll 1$ and $\|\boldsymbol{\gamma}\|\ll \sqrt{L}$. The estimation of the last sum follows the same arguments as above. Thus we obtain
\[
U_{rs}(L)\ll \frac{\ln{H}}{H}\cdot x^{4\varepsilon/5}\biggl\{\frac{x^{39/40}}{\Delta^{29/20}}\,L^{27/20} + \frac{x^{15/16}}{\Delta^{11/8}}\,L^{3/2}\biggr\}.
\]
Taking $L = P\cdot 2^{-\alpha}$ and summing over $\alpha$ we find:
\begin{multline*}
\sum\limits_{1\leqslant |h|\leqslant P}C(h)\sum\limits_{\substack{E<\kappa_{1},\kappa_{2}\leqslant K \\ (\kappa_{1},\kappa_{2})=1}}
e\biggl(\frac{h\varphi - h\nu\overline{\kappa}_{1}}{\kappa_{2}}\biggr)\ll \frac{x^{5\varepsilon/6}}{H}\biggl\{\frac{x^{39/40}}{\Delta^{29/20}}\,P^{27/20} + \frac{x^{15/16}}{\Delta^{11/8}}\,P^{3/2}\biggr\}\ll \\
\ll x^{\varepsilon}\biggl\{\frac{x^{39/40}}{\Delta^{29/20}}\,H^{7/20} + \frac{x^{15/16}}{\Delta^{11/8}}\,H^{1/2}\biggr\}.
\end{multline*}
Hence,
\[
U_{rs}^{(2)}\ll x^{\varepsilon}\biggl\{\frac{x^{39/40}}{\Delta^{29/20}}\,H^{7/20} + \frac{x^{15/16}}{\Delta^{11/8}}\,H^{1/2} + \frac{x}{\Delta^{2}H}\biggr\}.
\]
Summing this inequality over all pairs $r,s$ and using (\ref{lab_23}), we obtain:
\begin{multline*}
V(\varphi;\Delta,\xi_{1},\xi_{2};m_{1}, m_{2}) = U^{(1)}+U^{(2)} \ll \\
\ll x^{\varepsilon}\biggl\{\frac{x^{39/40}}{\Delta^{29/20}}\,H^{7/20} + \frac{x^{15/16}}{\Delta^{11/8}}\,H^{1/2}
+\frac{x}{\Delta^{2}H}\biggr\}\ll \\
\ll x^{1+\varepsilon}\biggl\{\frac{x^{-1/40}}{\Delta^{29/20}}\,H^{7/20} + \frac{x^{-1/16}}{\Delta^{11/8}}\,H^{1/2} + \frac{1}{\Delta^{2}H}\biggr\}.
\end{multline*}
The summation over $m_{i}\equiv\xi_{i}\pmod{\Delta}$, $i = 1,2$ will multiple the above bound to the factor $(\sqrt{x}/\Delta)^{2} = x/\Delta^{2}$:
\[
V(\varphi;\Delta,\xi_{1},\xi_{2}) \ll x^{2+\varepsilon}\biggl\{\frac{x^{-1/40}}{\Delta^{69/20}}\,H^{7/20} + \frac{x^{-1/16}}{\Delta^{27/8}}\,H^{1/2} + \frac{1}{\Delta^{4}H}\biggr\}.
\]
By Lemma 3,
\begin{multline*}
V(\varphi)\ll x^{2+\varepsilon}\sum\limits_{\Delta\leqslant D}\rho(\Delta)\biggl\{\frac{x^{-1/40}}{\Delta^{69/20}}\,H^{7/20} + \frac{x^{-1/16}}{\Delta^{27/8}}\,H^{1/2} + \frac{1}{\Delta^{4}H}\biggr\} \ll\, \\
\ll x^{2+\varepsilon}\biggl\{x^{-1/40}\,H^{7/20} + x^{-1/16}\,H^{1/2} +\frac{1}{H}\biggr\}.
\end{multline*}
Choosing $H = [x^{1/54}]$ we finally get
\[
V(\varphi)\ll x^{2+\varepsilon}\biggl(\frac{1}{H} + Hx^{-1/16}\biggr)\ll x^{2-1/54+\varepsilon}.
\]
Theorem is proved.
\vspace{0.5cm}

\textbf{\S 6. Multiple integral $\boldsymbol{\mathfrak{S}(\varepsilon)}$.}
\vspace{0.5cm}

In this section, we consider the integral $\mathfrak{S}(\varepsilon)$ defined in \S 4, that is,
\[
\mathfrak{S}(\varepsilon) =
\int_{0}^{1}du_{1}\int_{0}^{1}du_{2}\int_{\varepsilon}^{1+u_{1}}dv_{1}\int_{\varepsilon}^{1+u_{2}}\varphi(u_{1},u_{2},v_{1},v_{2})dv_{2},
\]
where we put
\begin{multline*}
\varphi(u_{1},u_{2},v_{1},v_{2}) = \frac{\max{(0,g - f)}}{v_{1}v_{2}},\quad g =\min{(g_{1},g_{2},g_{3})},\quad f = \max{(f_{1},f_{2},f_{3},f_{4})}, \\
f_{1} = u_{1}^{2},\quad f_{2} = u_{2}^{2},\quad f_{3} = (u_{1}-v_{1})^{2},\quad f_{4} = (u_{2}-v_{2})^{2}, \\
g_{1} = 1,\quad g_{2} = (u_{1}+v_{1})^{2},\quad g_{3} = (u_{2}+v_{2})^{2}
\end{multline*}
(in \S 3, we have already used the notations $f_{r}, g_{s}$ in other sense, but this will not lead us to misunderstandings).
Denote by $\omega$ the set of points $(u_{1},u_{2},v_{1},v_{2})$ satisfying the conditions $0\leqslant u_{i}\leqslant 1$, $0< v_{i}\leqslant 1+u_{i}$,  $i=1,2$.
Let us split $\omega$ to the subsets $\omega_{r,s}$ defined by the conditions $f = f_{r}$, $g = g_{s}$, $f_{r}<g_{s}$. Finally, denote by $\omega_{r,s}(\varepsilon)$
the set of points of $\omega_{r,s}$ with additional conditions $v_{i}>\varepsilon$, $i=1,2$.

Consequently, for any $0\leqslant \varepsilon < 1$ we get
\[
\mathfrak{S}(\varepsilon) = \sum\limits_{r,s}\int_{\omega_{r,s}(\varepsilon)}\varphi(u_{1},u_{2},v_{1},v_{2})\,du_{1}du_{2}dv_{1}dv_{2} =
\sum\limits_{r,s}\mathfrak{S}_{r,s}(\varepsilon).
\]
Thus, it is sufficient to show that
\begin{equation}\label{lab_24}
\mathfrak{S}_{r,s}(\varepsilon) = \mathfrak{S}_{r,s}(0) + O\biggl(\varepsilon\ln^{2}\frac{1}{\varepsilon}\biggr)
\end{equation}
for any pair $r,s$. We limit ourselves by the case $r=s=2$ since the arguments in all other cases follow the same lines.

One can check that $\omega_{2,2}$ is defined by the following set of inequalities:
\[
0\leqslant u_{1}\leqslant u_{2}\leqslant 1,\quad u_{2}-u_{1}\leqslant v_{1}\leqslant \min{(1-u_{1},u_{2}+u_{1})},\quad v_{1}-(u_{2}-u_{1})\leqslant v_{2}\leqslant 2u_{2}
\]
and, hence, $\omega_{2,2}$ splits into three subdomains $\varpi_{1}, \varpi_{2}$ and $\varpi_{3}$ defined by the systems
\begin{equation*}
\begin{cases}
0\leqslant u_{1}\leqslant u_{2}\leqslant \tfrac{1}{3},\\
u_{2}-u_{1}\leqslant v_{1}\leqslant u_{2}+u_{1},\\
v_{1}-(u_{2}-u_{1})\leqslant v_{2}\leqslant 2u_{2},
\end{cases}
\begin{cases}
\tfrac{1}{3} < u_{2} \leqslant 1,\\
0 \leqslant u_{1} \leqslant \tfrac{1}{2}(1-u_{2}),\\
u_{2}-u_{1}\leqslant v_{1}\leqslant u_{2}+u_{1},\\
v_{1}-(u_{2}-u_{1})\leqslant v_{2}\leqslant 2u_{2},
\end{cases}
\begin{cases}
\tfrac{1}{3} < u_{2} \leqslant 1,\\
\tfrac{1}{2}(1-u_{2}) < u_{1}\leqslant u_{2},\\
u_{2}-u_{1}\leqslant v_{1}\leqslant 1-u_{1},\\
v_{1}-(u_{2}-u_{1})\leqslant v_{2}\leqslant 2u_{2}.
\end{cases}
\end{equation*}
Denote by $\varpi_{j}(\varepsilon)$ ($j=1,2,3$) the set of points of $\varpi_{j}$ with additional condition $v_{i}>\varepsilon$ and denote by $J_{1}(\varepsilon)$ the integral
\begin{multline}\label{lab_25}
\int_{\varpi_{1}(\varepsilon)}\varphi\,\,du_{1}du_{2}dv_{1}dv_{2} = \\ = \int_{0}^{1/3}du_{2}\int_{0}^{u_{2}}du_{1}\int_{\max{(u_{2}-u_{1},\,\varepsilon)}}^{u_{2}+u_{1}}dv_{1}\int_{\max{(v_{1}-(u_{2}-u_{1}),\,\varepsilon)}}^{2u_{2}}\varphi\,du_{1}du_{2}dv_{1}dv_{2},\\
\varphi = \varphi(u_{1},u_{2},v_{1},v_{2}) = \frac{g_{2}-f_{2}}{v_{1}v_{2}} = \frac{(u_{1}+v_{1})^{2}-u_{2}^{2}}{v_{1}v_{2}}.
\end{multline}
From (\ref{lab_25}), it follows that $2u_{2}\geqslant \varepsilon$ for any point of $\varpi_{1}(\varepsilon)$. The contribution $j_{1}(\varepsilon)$ coming to $J_{1}(\varepsilon)$ from the subdomain with $\varepsilon/2\leqslant u_{2}\leqslant 2\varepsilon$, is less than
\[
\int_{\varepsilon/2}^{2\varepsilon}du_{2}\int_{0}^{2\varepsilon}du_{2}\int_{0}^{2\varepsilon}du_{1}\int_{\varepsilon}^{4\varepsilon}\frac{dv_{1}}{v_{1}}\int_{\varepsilon}^{4\varepsilon}\frac{(2\varepsilon + 4\varepsilon)^{2}}{v_{2}}\,dv_{2}\,\ll \varepsilon^{\,4}.
\]
Therefore,
\begin{multline*}
J_{1}(\varepsilon) = \int_{2\varepsilon}^{1/3}du_{2}\int_{0}^{u_{2}}du_{1}\int_{\max{(u_{2}-u_{1},\,\varepsilon)}}^{u_{2}+u_{1}}dv_{1}
\int_{\max{(v_{1}-(u_{2}-u_{1}),\,\varepsilon)}}^{2u_{2}}\varphi\,dv_{2} + O(\varepsilon^{\,4}) = \\
= \int_{2\varepsilon}^{1/3}du_{2}\int_{0}^{u_{2}-\varepsilon}du_{1}\int_{u_{2}-u_{1}}^{u_{2}+u_{1}}dv_{1}
\int_{\max{(v_{1}-(u_{2}-u_{1}),\,\varepsilon)}}^{2u_{2}}\varphi\,dv_{2} + \\
+ \int_{2\varepsilon}^{1/3}du_{2}\int_{u_{2}-\varepsilon}^{u_{2}}du_{1}\int_{\varepsilon}^{u_{2}+u_{1}}dv_{1}
\int_{\max{(v_{1}-(u_{2}-u_{1}),\,\varepsilon)}}^{2u_{2}}\varphi\,dv_{2} + O(\varepsilon^{4}) = \\
= J_{2}(\varepsilon) + J_{3}(\varepsilon) + O(\varepsilon^{\,4}).
\end{multline*}
Further,
\[
0\leqslant J_{3}(\varepsilon)\leqslant \int_{2\varepsilon}^{1/3}du_{2}\int_{u_{2}-\varepsilon}^{u_{2}}du_{1}\int_{\varepsilon}^{1}\frac{dv_{1}}{v_{1}}
\int_{\varepsilon}^{1}\frac{dv_{2}}{v_{2}} \ll \varepsilon\ln^{2}\frac{1}{\varepsilon}.
\]
Next, splitting the integral $J_{2}(\varepsilon)$ we get
\begin{multline*}
J_{2}(\varepsilon) = \int_{2\varepsilon}^{1/3}du_{2}\int_{0}^{u_{2}-\varepsilon}du_{1}\int_{u_{2}-u_{1}+\varepsilon}^{u_{2}+u_{1}}dv_{1}
\int_{v_{1}-(u_{2}-u_{1})}^{2u_{2}}\varphi\,dv_{2} + \\
+ \int_{2\varepsilon}^{1/3}du_{2}\int_{0}^{u_{2}-\varepsilon}du_{1}\int_{u_{2}-u_{1}}^{u_{2}-u_{1}+\varepsilon}dv_{1}
\int_{\varepsilon}^{2u_{2}}\varphi\,dv_{2} = J_{4}(\varepsilon) + J_{5}(\varepsilon).
\end{multline*}
Obviously,
\begin{multline*}
0\leqslant J_{5}(\varepsilon) \leqslant \int_{2\varepsilon}^{1/3}du_{2}\int_{0}^{u_{2}-\varepsilon}du_{1}\int_{u_{2}-u_{1}}^{u_{2}-u_{1}+\varepsilon}\frac{dv_{1}}{v_{1}}
\int_{\varepsilon}^{1}\frac{dv_{2}}{v_{2}} = \\ = \biggl(\ln{\frac{1}{\varepsilon}}\biggr)\int_{2\varepsilon}^{1/3}du_{2}\int_{0}^{u_{2}-\varepsilon}\ln{\biggl(\frac{u_{2}-u_{1}+\varepsilon}{u_{2}-u_{1}}\biggr)}\,du_{1}
= \biggl(\ln{\frac{1}{\varepsilon}}\biggr)\int_{2\varepsilon}^{1/3}du_{2}\int_{\varepsilon}^{u_{2}}\ln{\biggl(\frac{w+\varepsilon}{w}\biggr)}dw \leqslant \\ \leqslant \biggl(\ln{\frac{1}{\varepsilon}}\biggr)\int_{2\varepsilon}^{1/3}du_{2}\int_{\varepsilon}^{u_{2}}\frac{\varepsilon\,dw}{w}\leqslant
\varepsilon\biggl(\ln{\frac{1}{\varepsilon}}\biggr)^{2}\int_{2\varepsilon}^{1/3}du_{2}\ll \varepsilon\biggl(\ln{\frac{1}{\varepsilon}}\biggr)^{2}.
\end{multline*}
Now we replace the lower bound $u_{2}-u_{1}+\varepsilon$ for $v_{1}$ in $J_{4}(\varepsilon)$ to $u_{2}-u_{1}$. Thus we find
\begin{multline*}
J_{4}(\varepsilon) = \int_{2\varepsilon}^{1/3}du_{2}\int_{0}^{u_{2}-\varepsilon}du_{1}\int_{u_{2}-u_{1}}^{u_{2}+u_{1}}dv_{1}
\int_{v_{1}-(u_{2}-u_{1})}^{2u_{2}}\varphi\,dv_{2} - \\
- \int_{2\varepsilon}^{1/3}du_{2}\int_{0}^{u_{2}-\varepsilon}du_{1}\int_{u_{2}-u_{1}}^{u_{2}-u_{1}+\varepsilon}dv_{1}
\int_{v_{1}-(u_{2}-u_{1})}^{2u_{2}}\varphi\,dv_{2} = J_{6}(\varepsilon) - J_{7}(\varepsilon).
\end{multline*}
Setting $w = u_{2}-u_{1}$, $v = v_{1}-(u_{2}-u_{1})$ we have
\begin{multline*}
0\leqslant J_{7}(\varepsilon) \leqslant \int_{2\varepsilon}^{1/3}du_{2}\int_{\varepsilon}^{u_{2}}dw\int_{w}^{w+\varepsilon}\frac{dv_{1}}{v_{1}}
\int_{v_{1}-w}^{2/3}\frac{dv_{2}}{v_{2}} \leqslant \\
\leqslant \int_{2\varepsilon}^{1/3}du_{2}\int_{\varepsilon}^{1/3}dw\int_{0}^{\varepsilon}\frac{dv}{v+w}
\int_{v}^{2/3+w}\frac{dv_{2}}{v_{2}} \leqslant \\
\leqslant
\int_{2\varepsilon}^{1/3}du_{2}\int_{\varepsilon}^{1/3}\frac{dw}{w}\int_{0}^{\varepsilon}dv\int_{v}^{1}\frac{dv_{2}}{v_{2}} \leqslant \int_{0}^{1}du_{2}\int_{\varepsilon}^{1}\frac{dw}{w}\int_{0}^{\varepsilon}\biggl(\ln{\frac{1}{v}}\biggr)dv \ll \varepsilon\biggl(\ln{\frac{1}{\varepsilon}}\biggr)^{2}.
\end{multline*}
Next we replace the upper bound $u_{2}-\varepsilon$ for $u_{1}$ in $J_{6}(\varepsilon)$ to $u_{2}$. The error does not exceed
\[
J_{8}(\varepsilon) = \int_{2\varepsilon}^{1/3}du_{2}\int_{u_{2}-\varepsilon}^{u_{2}}du_{1}\int_{u_{2}-u_{1}}^{u_{2}+u_{1}}dv_{1}\int_{v_{1}-(u_{2}-u_{1})}^{2u_{2}}\varphi\,dv_{2}.
\]
Setting $w = u_{2}-u_{1}$, $v=v_{1}-(u_{2}-u_{1})$ we find that
\begin{multline*}
J_{8}(\varepsilon) = \int_{2\varepsilon}^{1/3}du_{2}\int_{0}^{\varepsilon}dw\int_{w}^{2u_{2}-w}dv_{1}
\int_{v_{1}-w}^{2u_{2}}\varphi\,dv_{2} = \int_{2\varepsilon}^{1/3}du_{2}\int_{0}^{\varepsilon}dw\int_{0}^{2(u_{2}-w)}dv
\int_{v}^{2u_{2}}\varphi\,dv_{2}.
\end{multline*}
At the same time,
\[
\varphi = \frac{(u_{1}+v_{1})^{2}-u_{2}^{2}}{v_{1}v_{2}} = \frac{(u_{1}+v_{1}-u_{2})(u_{1}+v_{1}+u_{2})}{v_{1}v_{2}} = \frac{v(v+2u_{2})}{v_{2}(v+w)}
\leqslant \frac{4u_{2}v}{v_{2}(v+w)}\leqslant \frac{4u_{2}}{v_{2}}.
\]
Therefore,
\begin{multline*}
J_{8}(\varepsilon) \leqslant 4\int_{2\varepsilon}^{1/3}u_{2}du_{2}\int_{0}^{\varepsilon}dw\int_{0}^{2u_{2}}dv
\int_{v}^{2u_{2}}\frac{dv_{2}}{v_{2}} \leqslant 4\varepsilon\int_{2\varepsilon}^{1/3}u_{2}\,du_{2}\int_{0}^{2u_{2}}\ln{\frac{2u_{2}}{v}}dv\leqslant \\
\leqslant 8\varepsilon \int_{0}^{1/3}u_{2}^{2}du_{2}\int_{1}^{+\infty}\frac{\ln{t}}{t^{2}}\,dt \ll \varepsilon.
\end{multline*}
Hence we have $J_{6}(\varepsilon) = J_{9}(\varepsilon) + O(\varepsilon)$, where
\[
J_{9}(\varepsilon) = \int_{2\varepsilon}^{1/3}du_{2}\int_{0}^{u_{2}}du_{1}\int_{u_{2}-u_{1}}^{u_{2}+u_{1}}dv_{1}\int_{v_{1}-(u_{2}-u_{1})}^{2u_{2}}\varphi\,dv_{2}.
\]
Finally, the error coming to $J_{9}(\varepsilon)$ from the replacement the lower bound $\varepsilon$ for $u_{2}$ by zero is estimated as follows ($w = u_{2}-u_{1}$, $v = v_{1}- (u_{2}-u_{1})$):
\begin{multline*}
\int_{0}^{2\varepsilon}du_{2}\int_{0}^{u_{2}}du_{1}\int_{u_{2}-u_{1}}^{u_{2}+u_{1}}dv_{1}
\int_{v_{1}-(u_{2}-u_{1})}^{2u_{2}}\varphi\,dv_{2} = \\
= \int_{0}^{2\varepsilon}du_{2}\int_{0}^{u_{2}}dw\int_{w}^{2u_{2}-w}dv_{1}\int_{v_{1}-w}^{2u_{2}}\varphi\,dv_{2}  =
\int_{0}^{2\varepsilon}du_{2}\int_{0}^{u_{2}}dw\int_{0}^{2(u_{2}-w)}dv\int_{v}^{2u_{2}}\varphi\,dv_{2}  \leqslant \\
\leqslant \int_{0}^{2\varepsilon}du_{2}\int_{0}^{u_{2}}dw\int_{0}^{2u_{2}}dv\int_{v}^{2u_{2}}\frac{4u_{2}}{v_{2}}\,dv_{2} \leqslant
4\int_{0}^{2\varepsilon}u_{2}du_{2}\int_{0}^{u_{2}}dw\int_{0}^{2u_{2}}\ln{\frac{2u_{2}}{v}}\,dv  \leqslant \\
\leqslant 8\int_{0}^{2\varepsilon}u_{2}^{2}du_{2}\int_{0}^{u_{2}}dw\int_{1}^{+\infty}\frac{\ln{t}}{t^{2}}\,dt\ll \varepsilon^{4}.
\end{multline*}
Thus we find that
\[
J_{1}(\varepsilon) = J_{1}(0) + O\biggl(\varepsilon\ln^{2}{\frac{1}{\varepsilon}}\biggr).
\]
The integrals over the domains $\varpi_{2}(\varepsilon), \varpi_{3}(\varepsilon)$ are considered along the same lines. Hence, the relation (\ref{lab_24}) is checked for $r=s=2$.
All other pairs $r,s$ are handle similarly.
\vspace{0.5cm}

\textbf{\S 7.Acknowledgements.}
\vspace{0.5cm}

The author is grateful to F.~Battistoni, L.~Greni\'{e} and G.~Molteni for informing him about the precise value of the constant $\mathfrak{S}$, and to the referees for careful reading of the text and for valuable remarks.

\end{document}